\documentclass[a4paper]{article}
\usepackage[utf8]{inputenc}
\usepackage[margin=2.5cm]{geometry} % full-width
\usepackage{multirow}
\usepackage{multicol}
\usepackage{amsmath}
\usepackage{amsthm}

\usepackage{amssymb}
\usepackage{tikz}
\usetikzlibrary{positioning, shapes}
\usetikzlibrary{topaths}
\usetikzlibrary{shapes.geometric}
\newcounter{algo@row}
\usetikzlibrary{arrows.meta}
\newcounter{algo@rowindent}
\newcommand{\algofont}[1]{\textbf{#1}}% S1
\newcommand{\algonumbersize}[1]{\scriptsize{#1}}% S2
\newcommand{\algopreitem}[1][\arabic{algo@row}]{\texttt{\algonumbersize{#1}}}
\newcommand{\algoitemskip}{\hspace{\value{algo@rowindent}cc}}
\newcommand{\algonewnestedopen}[2]{
	\newcommand{#1}[1][]{%
		\ifthenelse{\equal{##1}{}}{\item}{\item[{\algopreitem[##1]}]}
		\algoitemskip\algofont{#2}%
		\addtocounter{algo@rowindent}{1}%
		\ignorespaces
	}
}
\newcommand{\algonewnestedaux}[2]{
	\newcommand{#1}[1][]{
		\addtocounter{algo@rowindent}{-1}
		\ifthenelse{\equal{##1}{}}{\item}{\item[{\algopreitem[##1]}]}
		\algoitemskip\algofont{#2}%
		\addtocounter{algo@rowindent}{+1}%
		\ignorespaces
	}
}
\newcommand{\algonewnestedclose}[2]{
	\newcommand{#1}[1][]{
		\addtocounter{algo@rowindent}{-1}
		\ifthenelse{\equal{##1}{}}{\item}{\item[{\algopreitem[##1]}]}
		\algoitemskip\algofont{#2}%
		\ignorespaces
	}
}
\newcommand{\algonewcommand}[2]{
	\newcommand{#1}[1][default]{
		\ifthenelse{\equal{##1}{default}}{\item}{\item[{\algopreitem[##1]}]}%
		\algoitemskip\algofont{#2}%
		\ignorespaces
	}%
}
\newcommand{\algonewkeyword}[2]{\newcommand{#1}{\algofont{#2}}}
\algonewcommand{\STATE}{\ignorespaces}
\algonewcommand{\INPUT}{Input: }
\algonewcommand{\pINPUT}{\phantom{Input: }}
\algonewcommand{\COMPUTE}{Compute: }
\algonewcommand{\OUTPUT}{Output: }
\algonewcommand{\pOUTPUT}{\phantom{Output: }}
\algonewnestedopen{\IF}{if }
\algonewnestedaux{\ELSEIF}{else if }
\algonewnestedaux{\ELSE}{else }
\algonewnestedclose{\ENDIF}{end if }
\algonewnestedopen{\FOR}{for }
\algonewnestedclose{\ENDFOR}{end for }
\algonewnestedopen{\WHILE}{while }
\algonewnestedclose{\ENDWHILE}{end while }
\algonewcommand{\BREAK}{break}%
\algonewkeyword{\For}{for }%
\algonewkeyword{\To}{to }%
\algonewkeyword{\Do}{do }%
\algonewkeyword{\If}{if }%
\algonewkeyword{\Then}{then }%
\algonewkeyword{\Else}{else }%
\algonewkeyword{\End}{end }%
\algonewkeyword{\AND}{and }%
\algonewkeyword{\True}{true }%
\algonewkeyword{\False}{false }%
\algonewkeyword{\Call}{call }%
\algonewkeyword{\irbleigs}{irbleigs }%
\algonewkeyword{\tridiag}{tridiag}%
\algonewkeyword{\reorth}{reorth}%

\newtheorem{remark}{Remark}
\newcommand{\R}{\mathbb{R}}
\newcommand{\wtU}{\widetilde{\bU}}
\newcommand{\wtV}{\widetilde{\bV}}
\newcommand{\wtu}{\widetilde{\bu}}

\newcommand{\wtSigma}{\widetilde{\Sigma}}
\newcommand{\wtsigma}{\widetilde{\sigma}}
\newcommand{\wtLambda}{\widetilde{\Lambda}}
\newcommand{\wtlambda}{\widetilde{\lambda}}

\usepackage[utf8]{inputenc}
\usepackage{hyperref}
\hypersetup{
	unicode,
	colorlinks,
	urlcolor=blue, 
	linkcolor=blue, 
	citecolor=blue,
	pdfproducer={LaTeX},
	pdfcreator={pdflatex}
}
\usepackage{url}
\usepackage{xr-hyper}
\usepackage{listings}

\usepackage{pythonhighlight}
\definecolor{darkspringgreen}{rgb}{0.09, 0.45, 0.27}
\usepackage{enumitem}
\setlist{nolistsep}
\usepackage{xcolor}
\definecolor{lightergray}{rgb}{0.9,0.9,0.9}

\usepackage{xspace}
\usepackage{algorithm}
\usepackage{algpseudocode}
\usepackage[most]{tcolorbox}
\tcbuselibrary{listings,breakable}
\definecolor{darkspringgreen}{rgb}{0.09, 0.45, 0.27}

\usepackage{tikz}
\usetikzlibrary{arrows,shapes,positioning,shadows,trees}

\tikzset{
  basic/.style  = {draw, 
                   text width=6cm, 
                   font=\sffamily, 
                   rectangle},
  root/.style   = {basic, 
                   rounded corners=3pt, 
                   thin, 
                   align=center, 
                   fill=blue!50,
                   minimum height=0.7cm},
  level 2/.style = {basic, 
                    rounded corners=3pt, 
                    thin,align=center, 
                    fill=blue!23,
                    text width=1.85cm,
                    minimum height=0.5cm},
  level 3/.style = {basic, 
                    thin, 
                    align=left, 
                    fill=blue!5, 
                    text width=1.75cm}
}

\usepackage{soul}
\sethlcolor{lightergray}

\DeclareMathOperator*{\argmin}{arg\,min}
\usepackage{scalerel}
\usetikzlibrary{svg.path}
 
\definecolor{orcidlogocol}{HTML}{A6CE39}
\tikzset{
  orcidlogo/.pic={
    \fill[orcidlogocol] svg{M256,128c0,70.7-57.3,128-128,128C57.3,256,0,198.7,0,128C0,57.3,57.3,0,128,0C198.7,0,256,57.3,256,128z};
    \fill[white] svg{M86.3,186.2H70.9V79.1h15.4v48.4V186.2z}
                 svg{M108.9,79.1h41.6c39.6,0,57,28.3,57,53.6c0,27.5-21.5,53.6-56.8,53.6h-41.8V79.1z M124.3,172.4h24.5c34.9,0,42.9-26.5,42.9-39.7c0-21.5-13.7-39.7-43.7-39.7h-23.7V172.4z}
                 svg{M88.7,56.8c0,5.5-4.5,10.1-10.1,10.1c-5.6,0-10.1-4.6-10.1-10.1c0-5.6,4.5-10.1,10.1-10.1C84.2,46.7,88.7,51.3,88.7,56.8z};
  }
}

\newcommand\orcidicon[1]{\href{https://orcid.org/#1}{\mbox{\scalerel*{
\begin{tikzpicture}[yscale=-1,transform shape]
\pic{orcidlogo};
\end{tikzpicture}
}{|}}}}

\newcommand{\bSigma}{{\boldsymbol{\Sigma}}}
\newcommand{\bA}{{\bf A}}
\newcommand{\bB}{{\bf B}}
\newcommand{\bC}{{\bf C}}

\newcommand{\bE}{{\bf E}}
\newcommand{\bF}{{\bf F}}
\newcommand{\bG}{{\bf G}}
\newcommand{\bH}{{\bf H}}
\newcommand{\bI}{{\bf I}}

\newcommand{\bL}{{\bf L}}
\newcommand{\bM}{{\bf M}}

\newcommand{\bP}{{\bf P}}
\newcommand{\bQ}{{\bf Q}}
\newcommand{\bR}{{\bf R}}
\newcommand{\bS}{{\bf S}}

\newcommand{\bU}{{\bf U}}
\newcommand{\bV}{{\bf V}}
\newcommand{\bW}{{\bf W}}
\newcommand{\bX}{{\bf X}}
\newcommand{\bY}{{\bf Y}}
\newcommand{\bZ}{{\bf Z}}

\newcommand{\bb}{{\bf b}}
\newcommand{\bc}{{\bf c}}
\newcommand{\bd}{{\bf d}}
\newcommand{\be}{{\bf e}}
\newcommand{\bff}{{\bf f}}
\newcommand{\bg}{{\bf g}}

\newcommand{\bn}{{\bf n}}

\newcommand{\br}{{\bf r}}
\newcommand{\bs}{{\bf s}}
\newcommand{\bt}{{\bf t}}
\newcommand{\bu}{{\bf u}}
\newcommand{\bv}{{\bf v}}
\newcommand{\bw}{{\bf w}}
\newcommand{\bx}{{\bf x}}
\newcommand{\by}{{\bf y}}
\newcommand{\bz}{{\bf z}}
\newcommand{\bzero}{{\bf0}}

\newcommand{\calC}{\mathcal{C}}
\newcommand{\calD}{\mathcal{D}}

\newcommand{\calF}{\mathcal{F}}
\newcommand{\calG}{\mathcal{G}}

\newcommand{\calK}{\mathcal{K}}

\newcommand{\calR}{\mathcal{R}}
\newcommand{\calS}{\mathcal{S}}

\usepackage{listings}
\usepackage{xcolor}

\definecolor{codegreen}{rgb}{0,0.6,0}
\definecolor{codegray}{rgb}{0.5,0.5,0.5}
\definecolor{codepurple}{rgb}{0.58,0,0.82}
\definecolor{backcolour}{rgb}{0.95,0.95,0.92}

\lstdefinestyle{mystyle}{
    commentstyle=\color{codegreen},
    keywordstyle=\color{magenta},
    numberstyle=\tiny\color{codegray},
    stringstyle=\color{codepurple},
    basicstyle=\ttfamily \small,
    breakatwhitespace=false,         
    breaklines=true,                 
    keepspaces=true,                 
    numbersep=5pt,                  
    showspaces=false,                
    showstringspaces=false,
    showtabs=false,                  
    tabsize=2
}

\title{TRIPs-Py: Techniques for Regularization of Inverse Problems\\in Python}
\author{Mirjeta Pasha \orcidicon{0000-0003-4249-2421} \thanks{Department of Mathematics, Tufts University, Medford, MA, USA and Laboratory for Decision and Information System, Massachusetts Institute of Technology, Cambridge, MA, USA. (mpasha@mit.edu)}\and Silvia Gazzola \orcidicon{0000-0001-9588-0896} \thanks{Department of Mathematical Sciences, University of Bath, Bath, UK. (sg968@bath.ac.uk)}\and Connor Sanderford \thanks{Department of Mathematical Sciences, Arizona State University, Tempe, AZ, USA. (csanderf@asu.edu)}\and Ugochukwu O. Ugwu \thanks{Department of Electrical and Computer Engineering, Tufts University, Medford, MA, USA.(Ugochukwu.Ugwu@tufts.edu)}}
\begin{document}
\maketitle
\begin{abstract}
Abstract: In this paper, we describe TRIPs-Py, a new Python package of linear discrete inverse problems solvers and test problems. The goal of the package is two-fold: 1) to provide tools for solving small and large-scale inverse problems, and 2) to introduce test problems arising from a wide range of applications. The solvers available in TRIPs-Py include direct regularization methods (such as truncated singular value decomposition and Tikhonov) and iterative regularization techniques (such as Krylov subspace methods and recent solvers for $\ell_p$-$\ell_q$ formulations, which enforce sparse or edge-preserving solutions and handle different noise types). All our solvers have built-in strategies to define the regularization parameter(s). Some of the test problems in TRIPs-Py arise from simulated image deblurring and computerized tomography, while other test problems model realistic problems in dynamic computerized tomography.
Numerical examples are included to illustrate the usage as well as the performance of the described methods on the provided test problems. To the best of our knowledge, TRIPs-Py is the first Python software package of this kind, which may serve both research and didactical purposes.  
\end{abstract}

\section{Introduction} \label{sec:introduction}

Inverse problems arise whenever one wants to recover information about a hidden quantity from measurements acquired via a physical process (forward problem). In this paper, we are interested in linear discrete inverse problems, which can be formulated as linear systems of equations or linear least squares problems of the form
\begin{equation}\label{eq: linearEq}
\bA \bx = \bb,\qquad \min_{\bx}\|\bA \bx -\bb\|_2, 
\end{equation}
where $\bA \in \R^{m\times n}$ is a suitable discretization of the forward operator, {$\bx \in \R^n$ represents a quantity of interest, and $\bb \in \R^m$ is the available data, which is typically corrupted by some unknown perturbations (noise) $\be$}, i.e.,
\begin{equation}\label{eq: min}
\bb = \bb_{\rm true} + \be = \bA\bx_{\rm true} + \be\,.
\end{equation}

To keep the following derivations simple, we assume $m\geq n$; extensions to the $m<n$ case are often straightforward. Many important applications, such as medical, seismic, and satellite imaging require the solution of inverse problems; see, for instance, \cite{boas2001imaging, miller2012environmental, bennett1996generalized}. {Test problems in TRIPs-Py model deblurring (or deconvolution problems) and computerized tomography problems; while many realistic instances of these problems can be generated using the provided synthetic data, also instances of dynamic tomography problems that use real data are included.} 
For this kind of problems, the matrix $\bA$ is typically ill-conditioned, with singular values that gradually decay and cluster at zero. {The implication is that the solution of \eqref{eq: linearEq} is very sensitive to the noise in $\bb$} and some regularization should be applied to recover a meaningful approximation of $\bx_{\rm true}$. All the regularization methods considered in this paper and available within TRIPs-Py compute a regularized solution $\bx_{\rm reg}$ by (approximately) solving the following optimization problem
\begin{equation}\label{eq: genRegPb}
% \bx_{\rm reg} = \argmin_{\bx \in \R^{n}}\|\widetilde{\bA}\bx-\bb\|_p^p + \lambda\|\Psi\bx\|_q^q,\quad p,\,q>0,
\bx_{\rm reg} = \argmin_{\bx \in \calD\subseteq \R^{n}}\calF(\bx) + \alpha\calR(\bx),\quad \alpha>0,
\end{equation}
where $\calF(\bx)$ is a \emph{fit-to-data term} that typically involves the matrix $\bA$ (or a modification thereof) and the data $\bb$, $\alpha$ is a \emph{regularization parameter}, $\calR$ is a \emph{regularization term}, and $\calD$ is a set of \emph{constraints}. Most of the methods in this paper take $\calF(\bx)= \|\bA\bx-\bb\|_p^p$, with $p>0$, and $\calR(\bx)=\|\Psi\bx\|_q^q$, with $\Psi\in\R^{k\times n}$ and $q>0$. The choice of $p$ is dictated by the distribution of the noise $\be$ (e.g., $p = 2$ for Gaussian white noise, $p = 1$ for impulse noise), while the choices of $q$ and of the regularization matrix $\Psi$ are dictated by prior information available on $\bx_{\rm true}$ (e.g., if the gradient of $\bx_{\rm true}$ is smooth, one typically takes $q=2$ and $\Psi=\nabla$, i.e., a discretization of the gradient operator; if $\bx_{\rm true}$ is sparse one takes $\Psi$ as the identity and $0<q\leq 1$). Approaches for solving \eqref{eq: genRegPb} depend on the particular choices of $p,\,q,\,\Psi$ and the features of the problem (e.g., if $\bA$ is small or large scale). A summary of the possible choices of the functionals $\calF$ and $\calR$ supported within TRIPs-Py is provided in Table \ref{table: solvers}, and more details are provided in Section \ref{sec: solvers}. Within the TRIPs-Py solvers, the set of constraints $\calD$ is taken to be the whole domain $\R^n$ for small-scale problems, and an appropriate linear subspace of $\R^n$ for large-scale problems; choices for the latter are detailed in Section \ref{sec: projM}.

\begin{figure}
    \centering
\begin{tikzpicture}[
    nd/.style={draw, double, thick,fill=cyan!10,text width=2.3cm,
      align=flush center, minimum height=0.9cm, inner sep=0.6em,
      rounded corners},
    % arw/.style={thick,shorten >=1mm,shorten <=1mm,-Stealth}]
    % every node/.style={draw,thick,fill=green!5,text width=2.3cm,
    %   align=flush center,minimum height=0.9cm, inner sep=0.7em, rounded corners},
    every edge/.append style={ultra thick,shorten >=1mm,
      shorten <=1mm,-Stealth,bend left}
    ]
  % \draw[help lines] grid(15,12);
  \node[nd] at (8, 10) (TRIPs-Py) {TRIPs-Py};
  \node[nd] at (12, 8) (Demos) {Demos};
  \node[nd] at (4.5, 8) (Trips) {Trips};
  \node[nd] at (1, 6) (Solvers) {Solvers};
  \node[nd] at (1, 4) (TableSolvers) {see Table 1};
     \node[nd] at (8, 4) (Other) {Other};
    \node[nd] at (8, 2.5) (Parameter selection) {Reg. param};
\node[nd] at (4.5, 6) (Test problems) {Test problems};
\node[nd] at (4.5, 4) (Deblurring1D) {Deblurring 1D};
\node[nd] at (4.5, 2.5) (Deblurring) {Deblurring};
\node[nd] at (4.5, 1) (Tomography) {Tomography};
\draw[local bounding box=rect, dash dot dot, ultra thick]  (6.45,4.65)  rectangle (9.5,1.85);
\draw[local bounding box=rect, dash dot dot, ultra thick]  (2.85,4.65)  rectangle (6,0.35);
\node[nd,text width=1cm] at (7.1, 1) (DP) {DP};
\node[nd, text width=1cm] at (9, 1) (GCV) {GCV};
\node[nd] at (8, 6) (Utilities) {Utilities};
\node[nd] at (11, 6) (Data) {Data};
\node[nd] at (14, 6) (TableDemos) {see Table 3};
\draw[->,shorten >=1mm,
      shorten <=1mm,-Stealth]
    (4.5, 5.5) -- (4.5, 4.5);
\draw[->,shorten >=1mm,
      shorten <=1mm,-Stealth]
    (8, 2) -- (7, 1.5);
\draw[->,shorten >=1mm,
      shorten <=1mm,-Stealth]
    (8, 2) -- (9, 1.5);
\draw[->,shorten >=1mm,
      shorten <=1mm,-Stealth]
    (8, 5.5) -- (8, 4.5);
\draw[->,shorten >=1mm,
      shorten <=1mm,-Stealth]
    (12, 7.5) -- (11, 6.5);
    \draw[->,shorten >=1mm,
      shorten <=1mm,-Stealth]
    (12, 7.5) -- (13.5, 6.5);
\draw[->,shorten >=1mm,
      shorten <=1mm,-Stealth]
    (4.5, 7.5) -- (1.5, 6.5);
    \draw[->,shorten >=1mm,
      shorten <=1mm,-Stealth]
    (4.5, 7.5) -- (4.5, 6.5);
      \draw[->,shorten >=1mm,
      shorten <=1mm,-Stealth]
    (4.5, 7.5) -- (7.5, 6.5);
  \draw[->,shorten >=1mm,
      shorten <=1mm,-Stealth]
    (8, 9.5) -- (4.3, 8.5);
     \draw[->,shorten >=1mm,
      shorten <=1mm,-Stealth]
    (8, 9.5) -- (12, 8.5);
    \draw[->,shorten >=1mm,
      shorten <=1mm,-Stealth]
    (1, 5.5) -- (1, 4.5);
\end{tikzpicture}
    \caption{Overview of the TRIPs-Py's structure and contents. Most of the files available in the `Utilities' directory are auxiliary functions that can be used by the TRIPs-Py solvers, such as functions to set the regularization operators, or to display data and reconstructions.}
    \label{fig: TRIPs-Py_struct}
\end{figure}
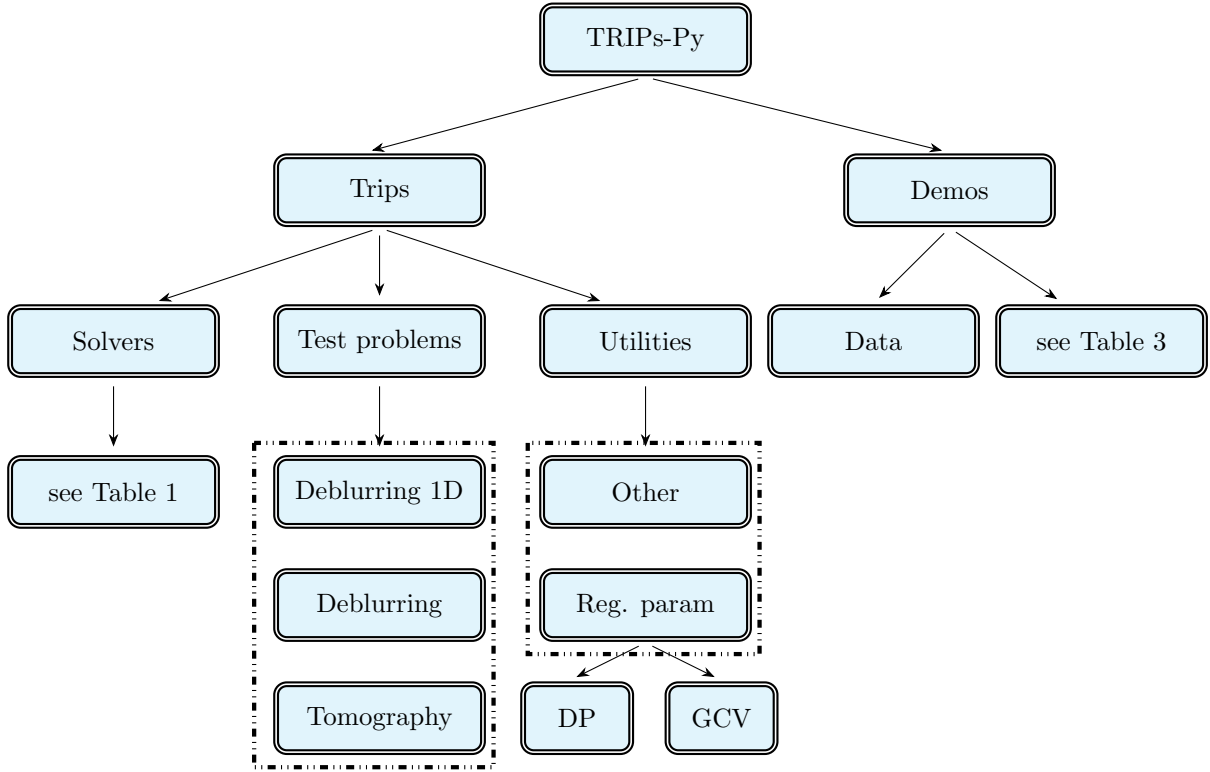

TRIPs-Py serves two interrelated purposes: 1) to provide model implementations of solvers for a variety of both small and large-scale linear inverse problems (including some recent methods that are not available in other packages), and 2) to provide a range of test problems (dealing with both synthetic and real data) for the users to test the TRIPs-Py solvers or, possibly, their own solvers. 
TRIPs-Py is an open source package that is available through GitHub at \url{https://github.com/trips-py/trips-py}, where the users can also find installation instructions and requirements. Figure \ref{fig: TRIPs-Py_struct} shows an overview of the TRIPs-Py structure and its modules.

When designing TRIPs-Py we aimed to create Python software that is easy to use, with calls to all the solvers that are very basic and similar. More precisely, default values are provided for all the options and parameters needed by the solvers (including automatic strategies to choose the regularization parameter in \eqref{eq: genRegPb} or the stopping iteration for iterative solvers); however, experienced users can easily set such parameters by themselves. The test problems in TRIPs-Py include 1D and 2D deconvolution problems (deblurring), and a variety of computed tomography problems that employ both synthetic and real data. All the test problems are generated with very similar instructions, and default values for the test problem generators are provided for users {who are} not familiar with the associated applications; as for the solvers, experienced users can easily set such problem-specific parameters by themselves.

Although other packages are already available in MATLAB and Python for solving, and numerically experimenting with, inverse problems (and TRIPs-Py shares some features with them), {some specifics} make TRIPs-Py unique. For instance, solely the solvers and test problems for small-scale inverse problems in TRIPs-Py are closely modeled on the ones available in \textsc{Regularization Tools} \cite{hansenRT}, a pioneering MATLAB package that has proven to be very popular with the wider community of researchers in algorithms for linear inverse problems since the Nineties. TRIPs-Py shares many design objectives, solvers and test problems with \textsc{IR Tools} \cite{gazzola2019ir}, a recent MATLAB package of iterative regularization methods and test problems for large-scale linear inverse problems. However, when compared to \textsc{IR Tools}, 
%regarding solvers, 
TRIPs-Py also features the first publicly available implementation of some recent methods for $\ell_p$-$\ell_q$ regularization and some test problems in TRIPs-Py also employ real data. Other popular MATLAB toolboxes such as \textsc{Restore Tools} \cite{NaPaPe04}, \textsc{AIR Tools II} \cite{hansenAIRII}) and \textsc{TIGRE} \cite{TIGRE} focus only on specific applications (the former is for image deblurring, the latter for computerized tomography), while TRIPs-Py encompasses a range of small and large-scale linear inverse problems. 
{Nowadays}, as Python is very popular with researchers both inside and outside of academia, and students at universities are increasingly exposed to Python as a programming language in their courses, we envision that TRIPs-Py could be useful for the community of users of these MATLAB toolboxes that may have to switch to Python for collaborating on research projects and for didactical purposes. Since many TRIPs-Py functionalities are similar to the ones underlying \textsc{Regularization Tools} and \textsc{IR Tools}, we are confident that such users will find transitioning to Python through TRIPs-Py natural. Even in the setting of the many powerful and popular Python packages for the solution of inverse problems, we believe that TRIPs-Py provides some valuable additions, for a variety of reasons. First, as mentioned above, TRIPs-Py's test problems model diverse applications, while many Python packages for inverse problems are focused on a particular application. This differentiates TRIPs-Py from recent packages like \textsc{CIL} \cite{CILI}, \textsc{ASTRA} \cite{van2015astra}, and the Python version of \textsc{TIGRE}, which all target computerized tomography. Second, TRIPs-Py's solvers for large-scale problems are based on standard or generalized Krylov subspace methods, with some of the former and the latter not being available anywhere else {or} being used to solve non-convex non-smooth instances of problem \eqref{eq: genRegPb}. This {distinguishes} TRIPs-Py from other general-purpose libraries for inverse problems like ODL \cite{ODL}, whose majority of solvers are based on optimization methods such as proximal gradient algorithms, primal-dual hybrid gradient algorithms, and ADMM. We anticipate that the users of these Python packages {would consider using some of the TRIPs-Py solvers to tackle their applications, and to compare available} or newly developed solvers with the ones in TRIPs-Py.

The remaining part of this paper is organized as follows: Section \ref{sec: solvers} gives an overview of the solvers available in TRIPs-Py, starting with filtering methods for small-scale problems and including many projection methods for large-scale, smooth and non-smooth regularized problems. Section \ref{sec: TestPbs} gives an overview of the test problems available in TRIPs-Py, together with some illustrations of the usage of solvers and test problems. Conclusions and future directions are discussed in Section \ref{sec: conclusions}.

\section{Overview of the TRIPs-Py solvers}\label{sec: solvers}

We start this section by discussing the TRIPs-Py solvers for regularization methods expressed in the 2-norm, i.e., to solve problem \eqref{eq: genRegPb} with $\calF(\bx)=\|\bA\bx-\bb\|_2^2$ and $\calR(\bx)=\|\Psi\bx\|_2^2$. This is the main part of the section, {and most of the TRIPs-Py solvers are tailored to this case}. The last part of this section describes a solver that can be employed for regularization methods expressed in the $\ell_p$-$\ell_q$ norm, i.e., to solve problem \eqref{eq: genRegPb} with $\calF(\bx)=\|\bA\bx-\bb\|_p^p$ and $\calR(\bx)=\|\Psi\bx\|_q^q$, $p,q>0$. An overview of the solvers available within TRIPs-Py is given in Table \ref{table: solvers}. We will not give too many details about the methods in this section, but rather provide extensive references for the reader. All the solvers in TRIPs-Py can be called in a consistent way that, at the very least, 
%(i.e., passing minimal amount of information) 
should include the forward operator (which could be a matrix or a function that acts on vectors) and the measured data. Additional optional inputs, such as the exact solution $\bx_{\rm true}$ for synthetic test problems, the maximum number of iterations to be performed and information about the stopping criterion for iterative solvers, or other specific features about specific methods, can be included. Most of these inputs are otherwise assigned default values. All functions return the computed approximation of the solution of the inverse problem, together with a dictionary that contains additional information about the solver and typically depends on the additional optional inputs assigned when calling the function. 

\begin{table}
  \centering
\caption{\label{table: solvers} List of the solvers available in TRIPs-Py.}
  \renewcommand{\arraystretch}{1.2}
  \begin{tabular}{|p{2.2cm}|p{5.7cm}|p{3.5cm}|p{1.6cm}|p{0.8cm}|}
    \hline
    \multirow{2}{3.5cm}{\textbf{Solver}} & 
    \multirow{2}{5cm}{\textbf{Description}} &
    \multicolumn{2}{l|}{\textbf{Specifics of problem \eqref{eq: genRegPb}}} & 
    \multirow{2}{1cm}{\textbf{Ref.}}\\
    \cline{3-4}
    & & $\calF(\bx)$ & $\calR(\bx)$ & \\
    %\hhline{~--}
    \hline\hline
    \texttt{TSVD} & Truncated SVD & $\|\bA_h\bx-\bb\|_2^2$ & 0 & \cite{hansen2010discrete} \\
    \texttt{TGSVD} & Truncated GSVD & $\|(\bA\Psi_{\bA}^{\dagger})_h\bx-\bb\|_2^2$ & 0 & \cite{hansen1998rank} \\
    \texttt{Tikhonov} & Tikhonov regularization & $\|\bA\bx-\bb\|_2^2$ & $\|\Psi\bx\|_2^2$ & \cite{hansen2010discrete} \\ 
    \texttt{CGLS} & Conjugate Gradient Least Squares & $\|\bA\bx-\bb\|_2^2$ & 0 & \cite{saad2003iterative} \\ 
    \texttt{GMRES} & Generalized Minimal Residual & $\|\bA\bx-\bb\|_2^2$, $\bA\in\R^{n\times n}$ & 0 & \cite{saad2003iterative} \\
    \texttt{Hybrid\_LSQR} & Hybrid LSQR & $\|\bA\bx-\bb\|_2^2$ & $\|\bx\|_2^2$ & \cite{review2023} \\
    \texttt{Hybrid\_GMRES} & Hybrid GMRES & $\|\bA\bx-\bb\|_2^2$, $\bA\in\R^{n\times n}$ & $\|\bx\|_2^2$ & \cite{review2023} \\
    \texttt{GK\_Tikhonov} & Golub-Kahan-Tikhonov & $\|\bA\bx-\bb\|_2^2$ & $\|\bx\|_{2}^2$ & \cite{fenu2016gcv}  \\
    \texttt{A\_Tikhonov} & Arnoldi-Tikhonov & $\|\bA\bx-\bb\|_2^2$ & $\|\bx\|_{2}^2$ & \cite{lewis2009arnoldi}  \\
    \texttt{GKS} & Generalized Krylov Subspace (GKS) & $\|\bA\bx-\bb\|_2^2$ & $\|\Psi\bx\|_2^2$ & \cite{lampe2012large} \\
    \texttt{MMGKS} & Majorization minimization with GKS & $\|\bA\bx-\bb\|_p^p$ & $\|\Psi\bx\|_q^q$ & \cite{lanza2015generalized} \\
    \texttt{AnisoTV} & Anisotropic Total Variation & $\|\bA\bx-\bb\|_2^2$ & $\|\nabla\bx\|_1$ & \cite{pasha2021efficient} \\
    \texttt{IsoTV} & Isotropic Total Variation & $\|\bA\bx-\bb\|_2^2$ & $\text{TV}(\bx)$ & \cite{pasha2021efficient}  \\
    \texttt{GS} & Group Sparsity & $\|\bA\bx-\bb\|_2^2$ & $\|\Psi\bx\|_{2,1}$ & \cite{pasha2021efficient}\\
    \hline
  \end{tabular}
\end{table}

We first survey methods that are suited for small-scale problems, {followed by} methods that are suited for large-scale problems; both the cases $\Psi=\bI$ and $\Psi\neq\bI$ will be covered. 

\subsection{Direct methods for small-scale problems}
When problem \eqref{eq: linearEq} is small scale, solvers for problems \eqref{eq: genRegPb} expressed in the 2-norm with $\calR(\bx)=\|\bx\|_2^2$ typically rely on the {S}ingular Value Decomposition (SVD) of $\bA$, i.e., %assuming $m\geq n$, 
$$
\bA = \bU\bSigma \bV^T,
$$
where $\bU\in\R^{m\times m}$ and
$\bV\in\R^{n\times n}$ are the orthogonal matrices of the left and right singular vectors, respectively, and $\bSigma
%= {\rm diag}(\sigma_1,\ldots,\sigma_m)
\in\R^{m\times n}$ is the matrix formed by the diagonal matrix of the singular values 
$\sigma_1\geq\sigma_2\geq\ldots\geq\sigma_n \geq 0$ {on top, and an $(m-n)\times n$ matrix of zeros at the bottom}. The SVD of $\bA$ is also a useful tool to analyze the ill-posedness of problem \eqref{eq: linearEq}.  While the SVD exists for every matrix, algorithms for its computation have a computational cost of order $O(mn^2)$, and are therefore prohibitive for large-scale problems; see 
\cite[\S 8.6]{GVL12} for more details. 

When $\calR(\bx)=\|\Psi\bx\|_2^2$, with $\bI\neq\Psi\in\R^{k\times n}$, problems \eqref{eq: genRegPb} are naturally handled by considering the generalized singular value decomposition (GSVD) of the matrix pair ($\bA$, $\Psi$). 
%, instead of the SVD of $\bA$.
{Assume that $m\geq n \geq k$, that $\text{rank}(\Psi)=k$ and that the null spaces of $\bA$ and $\Psi$ intersect trivially. Then the GSVD of} $(\bA$,$\Psi)$ is given by
\begin{equation}\label{eq:gsvd}
\bA=\wtU\wtSigma \bY^{T},\qquad \Psi=\wtV \wtLambda \bY^{T},
\end{equation}
where $\wtU\in\R^{m\times n}$ and $\wtV\in\R^{k\times n}$ have orthonormal columns,
$\bY\in\R^{n\times n}$ is nonsingular; $\wtSigma\in\R^{n\times n}$ is the diagonal matrix with diagonal entries $0\leq\wtsigma_1\leq\dots\leq\wtsigma_n\leq 1$, and $\wtLambda\in\R^{k\times n}$ is the matrix formed by the diagonal matrix of the values $1\geq\wtlambda_1\geq\dots\geq\wtlambda_k\geq 0$ on the left and $k\times (n-k)$ matrix of zeros on the right. The diagonal entries of $\wtSigma$ and $\wtLambda$ are such that, for $1\leq i\leq k$, $\wtsigma_i^2+\wtlambda_i^2=1$; the quantities $\wtsigma_i/\wtlambda_i$ are commonly referred to as the {\it generalized singular values} of $(\bA,\Psi)$. 
Similarly to the SVD of $\bA$, the cost of computing the GSVD of $(\bA,\Psi)$ is prohibitive for large-scale problems, unless some structure of $\bA$ and $\Psi$ can be exploited; see 
\cite[\S 6.1.6 and \S 8.7.4]{GVL12} for properties of the
GSVD and its computation. 

Next, we describe three TRIPs-Py basic regularization methods that are directly based on the SVD of $\bA$ or the GSVD of $(\bA,\Psi)$. Such methods are specific instances of the general class of (G)SVD filtering methods, whose solutions $\bx_{\mu}$ can be expressed as 
\begin{equation}\label{eq: filterM}
\bx_{\mu} = \sum_{i=1}^n\phi_i(\mu)\frac{\bu_i^T\bb}{\sigma_i}\bv_i\quad\mbox{(for $\Psi=\bI$)}\quad\mbox{and}\quad
\bx_{\mu} = \sum_{i=1}^k\phi_i(\mu)\frac{\wtu_i^T\bb}{\wtsigma_i}\widetilde{\by}_i + \sum_{i=k+1}^n(\wtu_i^T\bb)\widetilde{\by}_i\quad\mbox{(for $\Psi\neq\bI$)}\,,
\end{equation}
%for $\Psi=\bI$ and $\Psi\neq\bI$, respectively. 
where $\widetilde{\by}_i$, $i=1,\dots,n$ are the columns of $\bY^{-T}$. 
The scalars $\phi_i(\mu)$, $0\leq\phi_i(\mu)\leq 1$,  appearing in the above sums are {called {\it filter factors}. 
%Explicit expression of $\phi_i(\mu)$ for the (G)SVD will be provided below. 
The functional expressions of $\phi_i(\mu)$ determine different filtering methods, which all depend on the parameter $\mu$}. 
% ; these are detailed below for the TRIPs-Py solvers. 
The basic principles underlying filtering methods can be understood referring to the so-called `discrete Picard condition', which offers insight into the relative behavior of the magnitude of $\bu_i^T\bb$ and $\sigma_i$ (for $\Psi=\bI$), and $\wtu_i^T\bb$ and $\wtsigma_i$ (for $\Psi\neq\bI$), which appear in the expression of the filtered solutions. Namely, assuming that the quantities {$|\bu_i^T\bb_{\rm true}|$ and $\sigma_i$, and $|\wtu_i^T\bb_{\rm true}|$ and $\wtsigma_i$, decay at the same rate, then the noise in $\bb$ dominates the unregularized solution for small $\sigma_i$'s, and $\wtsigma_i$'s}. Therefore, filtering methods successfully compute regularized solutions when the filter factors $\phi_i(\mu)$ are close to 0 for small $\sigma_i$'s and $\wtsigma_i$'s, and close to 1 for large $\sigma_i$'s and $\wtsigma_i$'s (these are the meaningful components of the solution that we wish to retain). Having too many filter factors close to 1 results in under-regularized solutions (with $\phi_i(\mu)=1$ for every $i$ resulting in the (unregularized) solution of problem \eqref{eq: linearEq}); %with minimal $\|\Psi\bx\|_2$); 
having too many filter factors close to 0 results in over-regularized solutions. The second sum in the second equality in equation \eqref{eq: filterM} expresses the components of $\bx_{\mu}$ in the null space of $\Psi$, which are unaffected by regularization. Many strategies for choosing the regularization parameter $\mu$ are available in the literature, and Section \ref{sect:RegP} provides more details about the ones implemented in TRIPs-Py. 
In the following, we specify how the three filtering methods in TRIPs-Py relate to the general framework \eqref{eq: filterM}, specifying an expression for the filter factors $\phi_i(\mu)$.

\paragraph{Truncanted SVD (TSVD)}
The truncated SVD (TSVD) method regularizes \eqref{eq: linearEq} by computing 
\begin{equation}\label{eq: TSVD1}
\bx_h=\sum_{i=1}^{h} \frac{\bu_i^T\bb}{\sigma_i} \bv_i\,, \quad h\in\{1,2,\dots,n\}\,.
\end{equation}
TSVD is clearly a filtering method, with filter factors $\phi_i(h)=1$ for $1\leq i\leq h$ and $\phi_i(h)=0$ otherwise. Therefore, the integer $h$ plays the role of the regularization parameter, which can be set by applying the Generalized Cross-Validation (GCV) method or the discrepancy principle (if an estimate of the noise magnitude $\|\be\|_2$ is provided); see Section \ref{sect:RegP} for more details. {Let $\bA_h$ denote the best rank-$h$} approximation of $\bA$ in the 2-norm, i.e., 
\[
\bA_h = \bU_h\Sigma_h\bV_h^T,
\]
where $\bU_h$ and $\bV_h$ are obtained by taking the first $h$ columns of $\bU$ and $\bV$, respectively, and $\Sigma_h$ is the diagonal matrix with the first $h$ singular values on its main diagonal. Then $\bx_h$ can be regarded as the solution of the following variational problem
\[
\bx_h=\arg\min_{\bx\in\R^n}\|\bA_h\bx-\bb\|_2^2,
\]
which belongs to the framework \eqref{eq: genRegPb}, with $\calF(\bx)=\|\bA_h\bx-\bb\|_2^2$, $\calR(\bx)=0$, $\calD=\R^n$.

\paragraph{Tikhonov regularization}\label{subsec: standardformreg} Tikhonov regularization replaces \eqref{eq: linearEq} with the problem of computing 
\begin{equation}\label{eq: Tikhonov}
    \bx_{\alpha}=\arg\min_{\bx\in \R^{n}}\|\bA \bx - \bb\|_2^2 + \alpha \|\Psi \bx\|_2^2,\quad\mbox{where $\Psi\in\R^{k\times n}$}\,.
\end{equation}
When $\Psi = \bI$, problem \eqref{eq: Tikhonov} is said to be in the {\it standard form}, and when $\Psi \neq \bI$ it is said to be the {\it general form}.  Problem \eqref{eq: Tikhonov} clearly belongs to the framework \eqref{eq: genRegPb}, with $\calF(\bx)=\|\bA\bx-\bb\|_2^2$, $\calR(\bx)=\|\Psi\bx\|_2^2$, $\calD=\R^n$. Problem \eqref{eq: Tikhonov} can be equivalently expressed as a damped least squares problem, {with associated normal equations} 
\begin{equation}\label{eq: tikhclosesolL}
(\bA^T\bA + \alpha \Psi^T\Psi)\bx_{\alpha} = \bA^T\bb.
\end{equation}
If the null spaces of $\bA$ and $\Psi$ intersect trivially, the Tikhonov regularized solution $\bx_{\alpha}$ is unique. 
% \begin{equation}\label{eq: nanl}
% \mathcal{N}(\bA)\cap\mathcal{N}(\Psi)=\{0\}.
% \end{equation}
By plugging the SVD of $\bA$ (if $\Psi=\bI$) or the GSVD of $(\bA,\Psi)$ (if $\Psi\neq\bI$) into the above equation, one can see that $\bx_{\alpha}$ can be expressed in the framework of \eqref{eq: filterM} with 
\[
\phi_i(\alpha)=\frac{\sigma_i^2}{\sigma_i^2 + \alpha}\quad\mbox{(for $\Psi=\bI$)}\quad\mbox{and}\quad
\phi_i(\alpha)=\frac{(\wtsigma_i/\wtlambda_i)^2}{(\wtsigma_i/\wtlambda_i)^2 + \alpha}\quad\mbox{(for $\Psi\neq\bI$)}\,.
\]
From the above expression and looking at \eqref{eq: Tikhonov}, it is clear that more regularization is imposed for {larger} values of $\alpha$, {as more weight is put on the regularization term and more filter factors approach to $1$} (depending on the location of $\alpha$ within the range of the singular values of $\bA$ or generalized singular values of $(\bA,\Psi)$). Conversely, smaller values of $\alpha$ lead to under-regularized solutions. When $\Psi\neq\bI$ has a nontrivial null space, it is important to stress (again, from the second equation in \eqref{eq: filterM} or from \eqref{eq: Tikhonov}) that vectors in the null space of $\Psi$ are unaffected by regularization. The regularization parameter $\alpha$ can be set by applying the GCV method or the discrepancy principle (if an estimate of the noise magnitude $\|\be\|_2$ is provided); see Section \ref{sect:RegP} for more details. 
Finally, we remark that any Tikhonov regularized problem in general form can be equivalently transformed into a Tikhonov regularized problem in standard form. {The specific transformation depends on the properties of $\Psi$. Generically, one expresses} the equivalent standard form solution as
\begin{equation}\label{eq: stdfTikhonov}
\bx_{\alpha}=\Psi_{\bA}^{\dagger}\bz_{\alpha}+\bar{\bx}_{0},\quad\mbox{where}\quad    \bz_{\alpha}=\arg\min_{\bz\in \R^{k}}\|\bA\Psi_{\bA}^{\dagger} \bz - \bar{\bb}\|_2^2 + \alpha \|\bz\|_2^2.
\end{equation}
In the above expression, 
$\Psi^{\dagger}_{\bA} = (\bI - (\bA(\bI-\Psi^{\dagger}\Psi))^{\dagger}\bA)\Psi^{\dagger}$ denotes the $\bA$-weighted generalized pseudoinverse $\Psi^{\dagger}_{\bA}$ of the operator $\Psi$, $\bar{\bx}_0$ denotes the component of $\bx_{\alpha}$ in the null space of $\Psi$ and $\bar{\bb}=\bb-\bA\bar{\bx}_0$; see \cite{elden1982weighted} for more details.

\paragraph{Truncated GSVD (TGSVD)}
The truncated GSVD (TGSVD) method regularizes \eqref{eq: linearEq} by computing 
\begin{equation}\label{eq: TGSVD1}
\bx_h=\sum_{i=0}^{h-1} \frac{\wtu_{k-i}^T\bb}{\wtsigma_{k-i}} \widetilde{\by}_{k-i} + \sum_{i=k+1}^n(\wtu_i^T\bb)\widetilde{\by}_i\,,\quad h\in\{1,2,\dots,k\} \,.
\end{equation}
{TGSVD is a filtering method} that can be expressed as in \eqref{eq: filterM} (second equation), with filter factors $\phi_i(h)=1$ for $k-h+1\leq i\leq k$ and $\phi_i(h)=0$ otherwise. The TGSVD solution can be linked to both TSVD and Tikhonov regularization in general form, in that the TGSVD solution can {also be} expressed as
\begin{equation}\label{eq: TGSVD}
\bx_h=\Psi_{\bA}^{\dagger}\arg\min_{\bz\in\R^k}\|(\bA\Psi_{\bA}^{\dagger})_h\bz - \bb\|_2^2\,+\,\bar{\bx}_0\,,
% \quad\mbox{subject to}\quad\min_{\bx\in\R^n}\|\Psi\bx\|_2^2,    
\end{equation}
where $\Psi_{\bA}^{\dagger}$ and $\bar{\bx}_0$ are as in equation \eqref{eq: stdfTikhonov}, and $(\bA\Psi_{\bA}^{\dagger})_h$ is the optimal rank-$h$ approximation of $\bA\Psi_{\bA}^{\dagger}$ in the 2-norm, which can be either computed using the SVD of $\bA\Psi_{\bA}^{\dagger}$ or the GSVD of $(\bA,\Psi_{\bA})$. The variational problem \eqref{eq: TGSVD} belongs to the framework \eqref{eq: genRegPb}, as it is essentially a transformed least squares problem. As for TSVD, the truncation parameter $h$ plays the role of the regularization parameter, and similar strategies are used to automatically set it.

\subsection{Projection methods for large-scale problems}\label{sec: projM}

The methods described in the previous section can be properly applied only when it is feasible to compute some factorizations (e.g., the SVD) of the coefficient matrix $\bA$,  and possibly joint factorizations of the regularization matrix $\Psi$ and $\bA$ (e.g., the GSVD), making them suitable only for small scale problems or for problems where $\bA$, and possibly $\Psi$, have a special structure that can be exploited for storage and computations; see \cite{hansen2010discrete, hansen2006deblurring}. In general, when solving large-scale problems, only matrix-free methods, which do not require storage nor factorizations of $\bA$ but rather computations of matrix-vector products with $\bA$ and, often, {with $\bA^T$, are viable options}. %T(he majority of t)
This section is devoted to projection methods, which are iterative methods that {determine} a sequence of approximate regularized solutions of the original problem \eqref{eq: linearEq} in a sequence of subspaces of $\R^n$ of dimensions up to $d^\ast\ll n$. 
% (\uc{Are we assuming the dimensions of the subspaces are the same, and it is denoted by $d^*$?}). 
\begin{remark}\label{rem: projM}
All the projection methods available within TRIPs-Py compute a regularized solution by one of the following strategies:
\begin{enumerate}
    \item Applying a projection method directly to problems \eqref{eq: linearEq} and stopping before the noise corrupts the approximate solution (i.e., exploiting semiconvergence); see \cite[Chapter 6]{hansen2010discrete}.
    \item Applying a projection method to approximate the solution of a Tikhonov-like regularized problem. 
    \item Adopting a `hybrid' approach, whereby Tikhonov regularization is applied to a sequence of projected problems \cite{review2023}.
\end{enumerate}
As far as regularized problems in the 2-norm are concerned, all the approaches listed above can be employed. Moreover, when considering standard form Tikhonov regularization with a given regularization parameter, the second and third approaches are %mathematically 
equivalent. When considering regularization problems formulated in some $\ell_p$ norm, only the second approach will be considered within a computationally convenient 
%iteratively reweighting 
majorization-minimization strategy.
\end{remark}

In general, a projection method for computing an approximation $\bx_d$ of a solution of problems \eqref{eq: linearEq} is defined by the two conditions
\begin{equation}\label{eq: projM1}
%\[
\bx_d\in \calS_d,\quad \br_d := \bb - \bA\bx_d\perp \calC_d,
%\]
\end{equation}
where $\calS_d$ and $\calC_d$ are subspaces of $\R^n$ and $\R^m$ %, respectively, 
of dimension $d\ll \min\{m, n\}$, commonly referred to as the approximation subspace and the constraint subspace, respectively. If $\bS_d\in\R^{n\times d}$ and $\bC_d\in\R^{m\times d}$ are matrices whose columns span $\calS_d$ and $\calC_d$, respectively, then the above conditions can be equivalently expressed as
\begin{equation}\label{eq: projM2}
%\[
\mbox{find}\quad \bt_d\in\R^d\quad \mbox{such that}\quad
\bC_d^T(\bb - \bA\bS_d\bt_d)=\bzero\,,
%\]
\end{equation}
so that, since $d\ll n$, one can solve
\begin{equation}\label{eq: projpb}
\underbrace{\bC_d^T\bA\bS_d}_{=:\bG_d\in\R^{d\times d}}{\bt_d}=\underbrace{\bC_d^T\bb}_{=:\bg_d\in\R^d}
\end{equation}
by any direct method and form the approximate solution of \eqref{eq: linearEq} by taking $\bx_d=\bS_d\bt_d$. In practice, starting without loss of generality from the initial guess $\bx_0=\bzero$, a projection method for the solution of \eqref{eq: linearEq} is an iterative method that computes a sequence of approximate solutions $\{\bx_d\}_{d=1,2,\dots}$ satisfying conditions \eqref{eq: projM1} by building a sequence of nested subspaces
\[
% \calS_d=\calR ange(\bS_d)=\calR ange([\bS_{d-1}, \bs_d]),\quad
% \calC_d=\calR ange(\bC_d)=\calR ange([\bC_{d-1}, \bc_d]),\quad d=1, 2,\dots ,
\calS_d=\text{range}(\bS_d)=\text{range}([\bS_{d-1}, \bs_d]),\quad
\calC_d=\text{range}(\bC_d)=\text{range}([\bC_{d-1}, \bc_d]),\quad d=1, 2,\dots ,
\]
with $\bS_0, \bC_0$ being empty matrices; 
% . Without loss of generality, here and in the following we assume $\bx_0=\bzero$; 
see \cite{björck2014numerical} for more details. 

The majority of the solvers in TRIPs-Py are projection methods onto Krylov subspaces \cite{saad2003iterative}, i.e., both $\calS_d$ and $\calC_d$ are Krylov subspaces. In general, given $\bM\in \R^{n\times n}$ and $\bn\in \R^n$, a Krylov subspace is defined as
\begin{equation}\label{KrylovSp}
	\calK_d(\bM,\bn) = {\rm span}\{\bn, \bM \bn, \bM^2 \bn, \ldots , \bM^{(d-1)}\bn\}\,.
\end{equation}
Here and in the following, we assume that the dimension of $\calK_d(\bM,\bn)$ is $d$; appropriate checks are incorporated into TRIPs-Py to make sure that this assumption is satisfied. {Different Krylov methods differ in the choices of their approximation and constraint Krylov subspaces, and in their implementation}: when different Krylov subspace methods share the same $\calS_d$ and $\calC_d$, $d=1,2,\dots$, we say that they are mathematically equivalent. When Krylov methods are used to solve problems \eqref{eq: linearEq}, the quantities $\bM$ and $\bn$ appearing in \eqref{KrylovSp} are typically defined in terms of $\bA$ and $\bb$ (for square matrices), or $\bA^T\bA$ and $\bA^T\bb$, or $\bA\bA^T$ and $\bb$ (for rectangular matrices). In the following, we will refer to \eqref{KrylovSp} as standard Krylov subspace, to distinguish it from other Krylov-like subspaces that are generated by computing matrix-vector products with an iteration-dependent modification of the matrix $\bM$, and therefore cannot be expressed taking increasing powers of $\bM$ as in \eqref{KrylovSp}. 

\begin{table}
\renewcommand{\arraystretch}{1.5}
\caption{\label{table: decompositions} List of decompositions used by some TRIPs-Py solvers.}
\begin{tabular}{|p{38mm}|p{85mm}|p{15mm}|}\hline
\textbf{Decomposition} & \textbf{Description} &\textbf{Ref.} \\ \hline\hline
\texttt{arnoldi} & Arnoldi decomposition
& \cite{arnoldi1951principle} \\
\texttt{arnoldi\_update} & updated Arnoldi decomposition
& \cite{arnoldi1951principle} \\
\texttt{gsvd} & Generalized Singular Value Decomposition 
& \cite{GVL12}
\\
\texttt{golub\_kahan} & Golub-Kahan decomposition 
& \cite{GVL12}
\\
\texttt{golub\_kahan\_update} & updated Golub-Kahan decomposition 
& \cite{GVL12}
\\
\hline
\end{tabular}
\end{table}

\subsubsection{Methods based on standard Krylov subspaces}\label{sec: KS}

In this section, we briefly describe the Krylov subspace methods available within TRIPs-Py: namely, {Generalized Minimal Residual (GMRES)}(which can be applied when the coefficient matrix in \eqref{eq: linearEq} is square), {Least Squares QR (LSQR) and Conjugate Gradient Least Squares (CGLS)}, along with their available hybrid variants. We will emphasize how such methods relate to the framework in Remark \ref{rem: projM}.

\paragraph{Methods based on the Arnoldi algorithm: GMRES and its hybrid variant}
Given $\bA\in\R^{n\times n}$ and $\bb\in\R^n$, $d$ iterations of the Arnoldi algorithm initialized with $\bv_1 = \bb/\|\bb\|_2$ compute the partial factorization
\begin{equation}\label{eq: ArnoldiDecomp}
\bA \bV_d = \bV_{d+1} \bH_{d},
\end{equation}
where $\bV_{d+1} = [\bv_1, \bv_2,..., \bv_{d+1}]=[\bV_d, \bv_{d+1}] \in \R^{n \times (d+1)}$
has orthonormal columns that span the Krylov subspace $\mathcal{K}_{d+1}(\bA, \bb)$, and $\bH_{d} \in \R^{(d+1)\times d}$ is upper Hessenberg. TRIPs-Py provides two versions of the Arnoldi algorithm, both based on the modified Gram-Schmidt orthogonalization process (see Table \ref{table: decompositions}): one computes a given number of Arnoldi steps, while the other only computes one step, updating an already available partial Arnoldi decomposition. If the matrix $\bA$ is symmetric, the Arnoldi algorithm simplifies to the symmetric Lanczos algorithms \cite{lanczos1950iteration}, although only the former is provided in TRIPs-Py.  

{\emph{GMRES} is a Krylov method for the solution of the linear system appearing in \eqref{eq: linearEq}. At the $d$th iteration, GMRES} chooses $\calS_d=\calK_d(\bA,\bb)$, $\calC_d=\bA\calK_d(\bA,\bb)$ 
% \[
% \calS_d=\calK_d(\bA,\bb),\quad \calC_d=\bA\calK_d(\bA,\bb)
% \]
in \eqref{eq: projM1}. Because of this choice, {the associated residual $\br_d$ at the $d$th iteration is minimized in $\calK_d(\bA,\bb)$;} see \cite{saad2003iterative}. In practice, exploiting the Arnoldi factorization \eqref{eq: ArnoldiDecomp}, 
\begin{equation}% \label{eq: minBregKrylov}
\begin{array}{rcl}\label{eq: GMRES}
{\displaystyle\bx_d={\arg\min}_{\bx\in \mathcal{K}_d(\bA, \bb)}\|\bA\bx-\bb\|_{2}=                             
\displaystyle {\arg\min}_{\bt\in\mathbb{R}^{d}}\|\bA\bV_d\bt - \bb\|_{2}=
% \displaystyle\min_{\by\in\mathbb{R}^{d}}\|\bV_{d+1}\bH_{d+1,d}\by-\bb\|_{2}^{2}\\
% \overset{(b)}{=}
\displaystyle \arg\min_{\bt\in \mathbb{R}^{d}}\|\bH_{d}\bt-\|\bb\|_2 \be_{1}\|_{2},}
\end{array}
\end{equation}
where $\be_{1}$ denotes the first canonical basis vector of $\R^{d+1}$. 
Referring to the quantities introduced in equation \eqref{eq: projpb}, GMRES takes $\bG_d=\bH_d^T\bH_d$ and $\bg_d=\bH_d^T(\|\bb\|_2\be_1)$. 
It has been theoretically proven that GMRES is an iterative regularization method, so that the number of GMRES iterations serves as a regularization parameter; in other words, GMRES adopts the first strategy listed in Remark \ref{rem: projM}. If an estimate of the norm of the noise $\be$ affecting the data $\bb$ in \eqref{eq: min} is provided, 
%and assigned as optional input by the user, 
then TRIPs-Py automatically uses the discrepancy principle as stopping criterion for the GMRES iterations; see Section \ref{sect:RegP} for more details. Otherwise, the user may assign a maximum number of iterations to be performed, being aware that too many iterations may result in an under-regularized solution. From the first equality in equation \eqref{eq: GMRES}, we can see that GMRES can be expressed in the general framework of problem \eqref{eq: genRegPb}, with $\calF(\bx)=\|\bA\bx - \bb\|_2^2$, $\calR(\bx)=0$ and $\calD = \calK_{d^\ast}(\bA,\bb)$, where $d^\ast$ is the stopping iteration for GMRES. 

The \emph{hybrid GMRES} method applies some additional, iteration-dependent Tikhonov regularization in standard form to the projected problem appearing in \eqref{eq: GMRES}, following the third principle described in Remark \ref{rem: projM}. Namely, the $d$th iteration of the hybrid GMRES method computes the regularized solution
\begin{equation}
\begin{array}{rcl}\label{eq: hybridGMRES}
\displaystyle\bx_d=
% \arg\min_{\bx\in \mathcal{K}_d(\bA, \bb)}\|\bA\bx-\bb\|_{2}=                             
% \displaystyle \bV_d\arg\min_{\by\in\mathbb{R}^{d}}\|\bA\bV_d\by - \bb\|_{2}=
\displaystyle\bV_d\arg\min_{\bt\in \mathbb{R}^{d}}\|\bH_{d}\bt-\|\bb\|_2 \be_{1}\|_{2}^2 + \alpha_d\|\bt\|_2^2,
\end{array}
\end{equation}
where $\alpha_d$ is an iteration-dependent regularization parameter that can be automatically chosen applying the generalized cross validation method or the discrepancy principle (if an estimate of the noise magnitude $\|\be\|_2$ is provided) to the projected problem; we refer to Section \ref{sect:RegP} for more details. Assuming that a suitable value of $\alpha_d$ is chosen at each iteration, thanks to the additional regularization imposed in \eqref{eq: hybridGMRES} to the projected problem, hybrid GMRES is not affected by semiconvergence and setting a stopping criterion is less crucial than for GMRES; see \cite{review2023}. Therefore hybrid GMRES is stopped when a maximum number $d^\ast$ of iterations is performed. By exploiting decomposition \eqref{eq: ArnoldiDecomp} within \eqref{eq: hybridGMRES} and the definition of the matrix $\bV_{d+1}$, it can be shown that the hybrid GMRES method belongs to the general framework stated in \eqref{eq: genRegPb}, with $\calF(\bx)=\|\bA\bx-\bb\|_2^2$, $\calR(\bx)=\|\bx\|_2^2$, $\alpha=\alpha_{{d}^\ast}$ and $\calD=\calK_{{d}^\ast}(\bA,\bb)$. The \emph{Arnoldi-Tikhonov} method is mathematically equivalent to hybrid GMRES, the only difference between the two being that the former performs a given number of Arnoldi steps (which, if an estimate of $\|\be\|_2$ is available, is determined by the discrepancy principle), and regularizes the projected problem only at the last computed iteration.

\paragraph{Methods based on Golub-Kahan bidiagonalization (LSQR and its hybrid variant) and CGLS}

Given $\bA\in\R^{m\times n}$ and $\bb\in\R^m$, $d$ iterations of the Golub--Kahan bidiagonalization algorithm initialized with $\bv_1=\bA^T\bb/\|\bA^T\bb\|_2$ and $\bu_1=\bb/\|\bb\|_2$ compute the partial factorizations
\begin{equation}\label{eq: lanczosBidiag}
\bA\bV_{d}=\bU_{d+1}\bB_{d},\quad\bA^{T}\bU_{d+1}=\bV_{d+1}\bar{\bB}_{d+1}^{T},
\end{equation}
where  $\bV_{d+1} = [\bv_1, \bv_2,..., \bv_{d+1}]=[\bV_d, \bv_{d+1}] \in \R^{n \times (d+1)}$
has orthonormal columns that span the Krylov subspace $\mathcal{K}_{d+1}(\bA^T\bA, \bA^T\bb)$, $\bU_{d+1} = [\bu_1, \bu_2,..., \bu_{d+1}]=[\bU_d, \bu_{d+1}] \in \R^{m \times (d+1)}$
has orthonormal columns that span the Krylov subspace $\mathcal{K}_{d+1}(\bA\bA^T, \bb)$, $\bar{\bB}_{d+1} \in \R^{(d+1)\times (d+1)}$ is a lower bidiagonal matrix, and $\bB_{d}$ is obtained by taking the first $d$ columns of $\bar{\bB}_{d+1}$. As for the Arnoldi algorithm, 
TRIPs-Py provides two versions of the Golub--Kahan bidiagonalization algorithm (see Table \ref{table: decompositions}); both versions are implemented without reorthogonalization. 

\emph{LSQR} is a Krylov method for the solution of the problems appearing in \eqref{eq: linearEq} that, at the $d$th iteration, chooses $\calS_d=\calK_d(\bA^T\bA,\bA^T\bb)$, $\calC_d=\bA\calK_d(\bA^T\bA,\bA^T\bb)$ 
% \[
% \calS_d=\calK_d(\bA,\bb),\quad \calC_d=\bA\calK_d(\bA,\bb)
% \]
in \eqref{eq: projM1}. Because of this choice, at the $d$th iteration the residual $\br_d$ is minimized in $\calK_d(\bA^T\bA,\bA^T\bb)$; see \cite{saad2003iterative}. In practice, exploiting Golub-Kahan bigiagonalization \eqref{eq: lanczosBidiag}, 
\begin{equation}
\begin{array}{rcl}\label{eq: LSQR}
\displaystyle\bx_d=\arg\!\!\min_{\bx\in \mathcal{K}_d(\bA^T\bA, \bA^T\bb)}\!\!\|\bA\bx-\bb\|_{2}=                    \displaystyle \bV_d\arg\min_{\bt\in\mathbb{R}^{d}}\|\bA\bV_d\bt - \bb\|_{2}=         
\displaystyle\bV_d\arg\min_{\bt\in \mathbb{R}^{d}}\|\bB_{d}\bt-\|\bb\|_2 \be_{1}\|_{2}.
\end{array}
\end{equation}
Referring to the quantities introduced in equation \eqref{eq: projpb}, LSQR takes $\bG_d=\bB_d^T\bB_d$ and $\bg_d=\bB_d^T(\|\bb\|_2\be_1)$. It is well-known that LSQR is an iterative regularization method, so that the number of LSQR iterations serves as a regularization parameter; in other words, LSQR adopts the first strategy listed in Remark \ref{rem: projM}. Similarly to GMRES, the LSQR iterations are stopped when the discrepancy principle is satisfied (if an estimate of $\|\be\|_2$ is provided by the user) or if a maximum number of iterations is performed. From the first equality in equation \eqref{eq: LSQR}, we can see that LSQR can be expressed in the general framework of problem \eqref{eq: genRegPb}, with $\calF(\bx)=\|\bA\bx - \bb\|_2^2$, $\calR(\bx)=0$ and $\calD = \calK_{{d}^\ast}(\bA^T\bA,\bA^T\bb)$, where ${d}^\ast$ is the stopping iteration for LSQR. 

Similarly to hybrid GMRES, \emph{hybrid LSQR} applies some additional, iteration-dependent Tikhonov regularization in standard form to the projected problem appearing in \eqref{eq: LSQR}, following the third principle described in Remark \ref{rem: projM}. Namely, the $d$th iteration of the hybrid LSQR method computes the regularized solution
\begin{equation}
\begin{array}{rcl}\label{eq: hybridLSQR}
\displaystyle\bx_d=
% \arg\min_{\bx\in \mathcal{K}_d(\bA, \bb)}\|\bA\bx-\bb\|_{2}=                             
% \displaystyle \bV_d\arg\min_{\by\in\mathbb{R}^{d}}\|\bA\bV_d\by - \bb\|_{2}=
\displaystyle\bV_d\arg\min_{\bt\in \mathbb{R}^{d}}\|\bB_{d}\bt-\|\bA^T\bb\|_2 \be_{1}\|_{2}^2 + \alpha_d\|\bt\|_2^2,
\end{array}
\end{equation}
where $\alpha_d$ is an iteration-dependent regularization parameter that can be automatically chosen by applying the same approaches listed for hybrid GMRES. A proper choice of $\alpha_d$ mitigates the LSQR semiconvergence and the hybrid LSQR iterations are stopped when a maximum number ${d}^\ast$ of iterations is performed. Similarly to hybrid GMRES, hybrid LSQR belongs to the general framework stated in \eqref{eq: genRegPb}, with $\calF(\bx)=\|\bA\bx-\bb\|_2^2$, $\calR(\bx)=\|\bx\|_2^2$, $\alpha=\alpha_{{d}^\ast}$ and $\calD=\calK_{{d}^\ast}(\bA^T\bA,\bA^T\bb)$. Similarly to Arnoldi-Tikhonov, \emph{Golub-Kahan-Tikhonov} computes a fixed number of Golub-Kahan iterations and applies Tikhonov regularization only to the last computed projected problem.

\emph{CGLS} is a Krylov method for the solution of \eqref{eq: linearEq}, which is mathematically equivalent to LSQR. That is, at the $d$th iteration, $\calS_d$ and $\calC_d$ are as in LSQR, but the CGLS solution is computed using a three-term recurrence formula, rather than by (implicitly) computing the factorization \eqref{eq: lanczosBidiag} and solving the projected problem \eqref{eq: LSQR} as in LSQR. Both CGLS and LSQR can be naturally used to approximate the solution of the Tikhonov-regularized problem \eqref{eq: Tikhonov} (both in standard and general form, with a fixed regularization parameter $\alpha$), following the second strategy in Remark \ref{rem: projM}. This is achieved by applying such solvers to the equivalent formulation of \eqref{eq: Tikhonov} as a damped least squares problem, and stopping when we reach high accuracy in the associated normal equation residual. However, since CGLS does not compute the factorization \eqref{eq: lanczosBidiag}, it is not possible to state a hybrid version of CGLS as in \eqref{eq: hybridLSQR}.

\subsubsection{Methods based on generalized Krylov subspaces}\label{sec: GKS}
TRIPs-Py contains methods that rely on generalized Krylov subspaces, whose definition will be made clear in the next paragraph. Such methods include GKS, MMGKS, AnisoTV, IsoTV, and GS (see Table \ref{table: solvers} for a few more details). First we briefly describe the generalized Krylov subspace method (GKS) that approximates the solution of \eqref{eq: Tikhonov}, with a general regularization matrix $\Psi\in\R^{k\times n}$; see \cite{lampe2012large}. Then we describe the majorization-minimization (MM) technique that solves \eqref{eq: genRegPb} for a broad selection of $0<p,q \leq 2$ and $\Psi$, and we explain how GKS can be used to approximate the solution of the resulting reweighted problem: the resulting method is referred to as MMGKS. 
\paragraph{Generalized Krylov subspace (GKS) method}\label{sec: GKS}

The GKS method computes an approximation to \eqref{eq: Tikhonov} with a general $\Psi\in\R^{k\times n}$ by first determining an initial approximation subspace for the solution 
through $1\leq d\ll\min\{m,n\}$ steps of the Golub-Kahan bidiagonalization algorithm applied to $\bA$, with 
initial vector $\bb$. This gives a decomposition like the first one appearing in \eqref{eq: lanczosBidiag}. 
As far as $d$ is relatively small, it is inexpensive to compute the skinny QR factorizations
\begin{equation}\label{eq: skinnyQRl2l2}
\begin{array}{cccc}
\bA \bV_{d} = \bQ_{\bA} \bR_{\bA} \quad \text{with} &\quad \bQ_{\bA} \in \R^{m \times d}, & \bR_{\bA} \in \R^{d\times d},\\
\Psi\bV_{d} = \bQ_{\Psi} \bR_{\Psi} \quad \text{with} & \quad \bQ_{\Psi}\in \R^{r \times d}, & \bR_{\Psi} \in \R^{d \times d},
\end{array}
\end{equation}
where the matrices $\bQ_{\bA}$ and $\bQ_{\Psi}$ have orthonormal columns and the matrices $\bR_{\bA}$ and $\bR_{\Psi}$ are upper triangular. Taking $\bx_d=\bV_d\bt$ in \eqref{eq: Tikhonov} and using the factorizations in \eqref{eq: skinnyQRl2l2}, we obtain the %relatively small dimension 
$d$-dimensional linear system of equations
\begin{equation}\label{eq: normalEq_QR_l2l2}
(\bR^T_{\bA}\bR_{\bA} + \alpha \bR_{\Psi}^T\bR_{\Psi})\bt = \bR^T_{\bA}\bQ^T_{\bA}\bb,
\end{equation}
to be solved to get $\bt=\bt_d$. 
We then expand the approximation subspace for the solution of \eqref{eq: Tikhonov} by taking 

\begin{equation}\label{eq: enlargeMMGKS}%\label{eq: residual}
\bV_{d+1}=[\bV_{d},\bv_{\rm new}]\in\R^{n\times(d+1)},\quad\mbox{where}\quad\bv_{\rm new}=\widetilde{\br}_d/\|\widetilde{\br}_d\|_2,\quad\widetilde{\br}_d=\bA^T(\bA\bV_{d}\bt_d -\bb)+\alpha  \Psi^T\Psi\bV_{d}\bt_d,
\end{equation}
i.e., by adding the current normalised residual of the (full) normal equations. This concludes the first GKS iteration. 

Next, the QR factorization in \eqref{eq: skinnyQRl2l2} is updated for the matrices $\bA\bV_{d+1}$ and $\Psi\bV_{d+1}$ and the process in steps \eqref{eq: skinnyQRl2l2} - \eqref{eq: enlargeMMGKS} is repeated, expanding the 
% initial $\ell$-dimensional 
solution space during subsequent GKS iterations until a sufficiently accurate solution is reached. Since a suitable value of the regularization parameter $\alpha$ in \eqref{eq: Tikhonov} is usually not know and applying some parameter choice strategy to \eqref{eq: Tikhonov} is computationally unfeasible for large-scale problems, the GKS method allows an iteration-dependent choice of $\alpha$. Specifically, in TRIPs-Py one can adopt the discrepancy principle (if a good estimate of the noise magnitude $\|\be\|_2$ is provided by the user) or generalized cross validation to compute $
\alpha=\alpha_d$ for problem \eqref{eq: normalEq_QR_l2l2}; see Section \ref{sect:RegP} for more details. Since, in general,  $\alpha$ may vary at each iteration and $\Psi\neq\bI$, the range of $\bV_{d+1}$ (approximation subspace for the solution) is not a standard Krylov subspace anymore, and it is called generalized Krylov subspace. If $\alpha$ is fixed or $\Psi=\bI$, the GKS method is mathematically equivalent to LSQR or CGLS applied to the damped least squares formulation of the Tikhonov problem \eqref{eq: Tikhonov}.  
\paragraph{Majorization-minimization Generalized Krylov subspace (MMGKS) method} Let us consider problem \eqref{eq: genRegPb} with $\calF(\bx)=\|\bA\bx-\bb\|_p^p$ and $\calR(\bx) = \|\Psi\bx\|_q^q$. Due to the non-differentiability and possibly nonconvexity for certain choices of $p$ and $q$, it is common to replace the data fidelity and the regularization terms with differentiable approximations thereof, and solve the following problem
\begin{equation} \label{eq:Je}
\min_{\bx} \mathcal{J}_{\epsilon,\alpha}(\bx) = \min_{\bx} 
  \sum_{j = 1}^{m} \phi_{p, \epsilon}(( \bA \bx - \bb)_j) +  \alpha\sum_{j = 1}^{k} \phi_{q, \epsilon}(( \Psi\bx)_j)\,,
\end{equation}
where\begin{equation}\label{eq:smothedPhi}
 \phi_{r,\epsilon}(t)= \left(t^2+\epsilon^2\right)^{r/2}\mbox{~~with~~}
 \left\{ \begin{array}{ll}
\epsilon >0\mbox{~~for~~}0<r\leq 1,\\
\epsilon =0\mbox{~~for~~}r>1,
 \end{array}
 \right.
\end{equation}
for a small positive constant $\epsilon>0$. 
A well-known approach to solve \eqref{eq:Je} is majorization-minimization (MM) \cite{lange2016mm, rodriguez2008efficient}, which solves a sequence of regularized least squares problems. Namely, at the $(\ell+1)$th MM iteration, one computes 
\begin{equation} \label{eq:reg2}  \bx^{(\ell+1)} = \arg \min_{\bx\in\R^n} \|\bP_{\rm fid}^{(\ell)}\left(\bA \bx - \bb\right) \|_2^2 + \alpha\| \bP_{\rm reg}^{(\ell)}\Psi \bx \|_2^2 \, , 
\end{equation}
where $\bP_{s}^{(\ell)}$ ($s$ = fid, reg) denotes a diagonal weighting matrix, whose entries are determined from the current solution $\bx^{(\ell)}$, $\ell=0,1,\dots$. 

Classical methods for MM involve solving \eqref{eq:reg2} with a fixed regularization parameter by, e.g., applying CGLS, which results in time-consuming inner-outer iterative strategies. The MMGKS method implemented in TRIPs-Py simultaneously computes a new approximation $\bx^{(\ell+1)}$ and updates the weights $\bP_{s}^{(\ell+1)}$ ($s = \rm fid, \rm reg$) by projecting the current problem \eqref{eq:reg2} onto a generalized Krylov subspace, whereby an iteration-dependent suitable value for the regularization parameter $\alpha=\alpha_{\ell}$ can be computed (similarly to the method described in the previous paragraph) and the quantities $\bP_{\rm fid}^{(\ell)}\bA$, $\bP_{\rm fid}^{(\ell)}\bb$, and $\bP_{\rm reg}^{(\ell)}\Psi$ can be incorporated. In this way, for the MMKS methods, no inner-outer iteration scheme is needed and an MM iterate corresponds to a GKS iterate. Namely, the generalized Krylov subspace is expanded by 1 vector at each iteration, so that the $(\ell+1)$th approximate solution $\bx^{(\ell+1)}$ belongs to a generalized Krylov subspace of dimension $d=d_0+\ell$, where $d_0$ denotes the dimension of the initial approximation subspace for building the generalized Krylov subspace. We refer to  
\cite{huang2017majorization, lanza2015generalized}
%\cite{huang2017majorization, buccini2020modulus}
for further details. 

A variety of regularized formulations can be handled more or less straightforwardly by the MMGKS methods, provided that the parameters $p$ and $q$ (for the $\ell_p$ and $\ell_q$ (quasi) norms), the regularization matrix $\Psi$ and the weights $\bP_s^{(\ell)}$ ($s$ = fid, reg) are properly defined. Some of the solvers listed in Table \ref{table: solvers} are indeed drivers for MMGKS, whereby MMGKS is called with specific inputs, to conveniently handle particularly relevant regularization terms. Namely, 
\begin{itemize}
\item \texttt{AnisoTV} implements anisotropic total variation (TV) regularization. In this instance, $p=2$ and the matrix $\Psi\in\R^{k\times n}$ is a rescaled finite difference discretization of the first derivative operator in one or two spatial dimensions, or in three spatio-temporal dimensions; $k$ depends on the dimensionality of the problem. The weights are given by the diagonal matrix with entries 
\[
\left((\Psi\bx^{(\ell)})_j^2 + \epsilon^2\right)^{(q-2)/4},\quad\mbox{where}\quad 0<q\leq 1\,.
\]
% More details about this solvers are provided in \cite{huang2017majorization, lanza2015generalized}.
\item \texttt{IsoTV} implements isotropic TV regularization. In one spatial dimension this coincides with anisotropic TV, so that isotropic TV is more meaningful in two spatial dimensions or in three spatio-temporal dimensions. 
\item \texttt{GS} implements a regularizer enforcing group sparsity (possibly under transform), for solutions that are that are naturally partitioned in subsets; see~\cite{bach2012optimization}.
\end{itemize}
We refer to Section \ref{sec: DemoDynamic} for more details about (an)isotropic TV and group sparsity within MMGKS; see also \cite{pasha2021efficient}.

\subsection{Strategies to choose the regularization parameter}\label{sect:RegP}

In this section we discuss methods available within TRIPs-Py to choose the regularization parameter for the solvers described in the previous two subsections (and summarized in Table \ref{table: solvers}). Properly choosing regularization parameters is a pivotal task for any regularization method, as the success of the latter crucially depends on the former. In this section, when considering projection methods, it is  convenient to use the following compact notation for the regularized solution associated to the $d$th projected problem
\begin{equation}\label{eq: projregsol}
\bx_d=\bV_d\underbrace{\bF_{\alpha_d}^\dagger\bff_d}_{=:\bt_d},
\end{equation}
where $\bV_d$ is the basis for the $d$th approximation subspace, $\bF_{\alpha_d}^{\dagger}$ is the so-called $d$th regularized inverse for the projected problem (depending on the $d$th Tikhonov regularization parameter $\alpha_d$), and $\bff_d$ is the $d$th projected right-hand side vector. All quantities appearing in \eqref{eq: projregsol} are specific for each projection method. In particular, for projections methods that do not involve any additional Tikhonov regularization,  $\bx_d=\bV_d\bF_{0}^\dagger\bff_d$.

\paragraph{Discrepancy principle} 

Let us assume that an estimate 
$\delta$ for the norm of the noise $\be$ appearing in \eqref{eq: min} is available. Then the discrepancy principle computes the regularization parameter for the regularized solution of the generic regularization problem \eqref{eq: genRegPb} by imposing
\begin{equation}\label{errbd}
\calD(\bx_{\rm reg}):=\|\bA\bx_{\rm reg}-\bb\|_2=\eta\delta,
\end{equation}
where $\eta>1$ ($\eta\simeq 1$) is a safety factor. Since the specific expression and properties of the functional $\calD(\bx_{\rm reg})$ depend on the considered regularization method, TRIPs-Py includes different versions of \eqref{errbd}. In general, for the discrepancy principle to be satisfied, one should make the natural assumption that
\begin{equation}\label{eq: condDP1}
\|\bb_0\|_2\leq \|\be\|_2\leq\|\bb\|_2\,,  
\end{equation}
where $\bb_0$ denotes the orthogonal projection of $\bb$ onto the null space of $\bA\bA^T$; see \cite{ReichelAnzf}.

For T(G)SVD, $\calD(\bx_h)$ is a decreasing functional of the discrete truncation parameter $h$: therefore one chooses the truncation parameter $h^\ast$ such that 
\[
\calD(\bx_{h^\ast})\geq \eta\delta> \calD(\bx_{h^\ast+1})\,.
\]

Similarly, for purely iterative regularizing methods, $\calD(\bx_{\rm reg})$ is evaluated at discrete points (i.e., the number of iterations or, equivalently, the dimension of the projection subspace). For the $d$th projected problem, $\calD(\bx_d)$ is essentially the norm of the residual associated to the $d$th solution. For GMRES and LSQR, $\calD(\bx_d)$ is a decreasing functional of $d$: this is a consequence of the optimality properties mentioned in Section \ref{sec: KS}. Therefore, to satisfy the discrepancy principle, one stops at the iteration $d^\ast$ such that
\[
\calD(\bx_{d^\ast-1})\geq \eta\delta> \calD(\bx_{d^\ast})\,.
\]
Note that, for GMRES and LSQR, the discrepancy principle is computationally cheap to apply, as the functional $\calD(\bx_d)$ can be computed with respect to the projected coefficient matrix and right-hand side vector; see the last two equalities in equations \eqref{eq: GMRES} and \eqref{eq: LSQR}, respectively, for a justification.

When the discrepancy principle is applied to Tikhonov regularization (for either the full-dimensional problem or at each iteration of a projected problem), \eqref{errbd} amounts to a nonlinear equation to be solved with respect to the regularization parameter $\alpha$. Note that, with the change of variable $\hat{\alpha}=1/\alpha$ in \eqref{eq: Tikhonov}, \eqref{eq: normalEq_QR_l2l2}, \eqref{eq:reg2}, and $\hat{\alpha}_d=1/\alpha_d$ in \eqref{eq: hybridGMRES}, \eqref{eq: hybridLSQR}, $\calD(\bx_{\rm reg})$ is a decreasing convex functional of $\hat{\alpha}$ and $\hat{\alpha}_d$, respectively, so that Newton's method is guaranteed to converge when started on the left of the zero of $\calD(\bx_{\rm reg})-
\eta\delta$. For hybrid GMRES and hybrid LSQR, i.e., when considering the projected problems appearing in \eqref{eq: hybridGMRES} and \eqref{eq: hybridLSQR}, the discrepancies computed with respect to the full-dimensional and the projected problems coincide, i.e.,
\begin{equation}\label{eq: hybridDP}
\calD(\bx_{\rm reg})=\|\bb-\bA\bV_d\bF_{\alpha_d}^\dagger\bff_d\|_2=
\|\bff_d-\bF_{\alpha_d}^\dagger\bff_d\|_2\,;    
\end{equation}
see, again, the last two equalities in equations \eqref{eq: GMRES} and \eqref{eq: LSQR}. Therefore, in order for $\calD(\bx_{\rm reg})-
\eta\delta$ to have a zero, one should also assume that
\begin{equation}\label{eq: condDP2}
    \|(\bI-\bF_0^\dagger)\bff_d\|_2\leq\|\be\|_2\leq\|\bb\|_2\,,
\end{equation}
(it should probably be something like $(\bI - \bF\bF_0^\dagger)$) 
which is essentially condition \eqref{eq: condDP1} applied to the projected problems appearing in \eqref{eq: hybridGMRES} and \eqref{eq: hybridLSQR}. Since the quantity on the left of the first inequality is the norm of the GMRES or the LSQR residuals (recall that, with $\alpha_d=0$ hybrid GMRES or hybrid LSQR are equivalent to GMRES or LSQR, respectively), and since such residuals decrease with $d$, the above condition is satisfied when $d$ is sufficiently large (typically after only a few iterations); see also \cite{gazzola2020survey}. For methods based on generalized Krylov subspaces, the discrepancies computed with respect to the full-dimensional and the projected problems are different. Namely, 
\begin{equation}\label{eq: GKSDP}
\calD(\bx_{\rm reg})=\left(\|\bff_d- \bR_{\bA}\bF_{\alpha}^{\dagger}\bff_d\|_2^2 + \|(\bI - \bQ_{\bA}\bQ_{\bA}^T)\bb\|_2^2\right)^{1/2}
\end{equation}
when computed for the full-dimensional problem. The second term in the above sum is dropped when $\calD(\bx_{\rm reg})$ is computed with respect to the projected problem; see also \cite{BRXX}. In this setting, in order for $\calD(\bx_{\rm reg})-\eta\delta$ to have a zero, a condition similar to \eqref{eq: condDP2} should be satisfied; for the full-dimensional problem, the leftmost quantity in \eqref{eq: condDP2} should be replaced by $(\|(\bI-\bF_0^\dagger)\bff_d\|_2+\|(\bI - \bQ_{\bA}\bQ_{\bA}^T)\bb\|_2^2)^{1/2}$. By default, for the GKS-based solvers, TRIPs-Py applies the discrepancy principle with respect to the full-dimensional problem. 

\paragraph{Generalized Cross Validation (GCV)}

The GCV criterion prescribes to take the regularization parameter that minimizes the functional 
\begin{equation}\label{eq: GCV}
\calG(\bx_{\rm reg})=\frac{\|\bA\bx_{\rm reg}-\bb\|_2^2}{\left(\text{trace}(\bI - \bA\bA_{\rm reg}^{\dagger})\right)^2}\,.
\end{equation}
Such procedure is derived from statistical 
%considerations, 
techniques, starting from the principle that a good regularized solution $\bx_{\rm reg}$ (defined by a good regularization parameter) should be able to predict the %missing 
exact data %points in 
$\bb_{\rm true}$ as well as possible. In equation \eqref{eq: GCV}, $\bA_{\rm reg}^{\dagger}$ is the regularized inverse of $\bA$ (specific for each regularization method), i.e., a matrix such that $\bx_{\rm reg}=\bA_{\rm reg}^{\dagger}\bb$, and the quantity $\bA\bA_{\rm reg}^{\dagger}$ is often referred to as `influence matrix'. Since GCV does not require any information about the magnitude of the noise affecting the data $\bb$, it is the default regularization parameter choice method for the TRIPs-Py solvers that involve TSVD or Tikhonov regularization (for either the full-dimensional or the projected problem).

For (G)SVD spectral filtering methods, the functional $\calG(\bx_{\rm reg})$ can be conveniently expressed with respect to the filter factors and quantities appearing in the SVD of $\bA$ or the GSVD of $(\bA,\Psi)$. In particular, for T(G)SVD, the functional $\calG(\bx_h)$ is evaluated at discrete points and its denominator simplifies to $(m-h)^2$, with $h=1,\dots,n$ for TSVD and $h=1,\dots,k$ for GSVD (see equations \eqref{eq: TSVD1} and \eqref{eq: TGSVD1}, respectively); we refer to \cite[\S 5.4]{hansen2010discrete} for further details. 

When GCV is applied to the hybrid-projection methods based on standard Krylov subspace methods (Section \ref{sec: KS}) and to the methods based on generalized Krylov subspaces (Section \ref{sec: GKS}), the values of $\calG(\bx_{\rm reg})$ may depend on whether the regularization parameter is selected for the projected problem only, or by linking the projected problem to the corresponding full-dimensional regularized problem. 
%%%
% Let us consider the $d$th projected problem, whose regularized solution can be compactly expressed as 
% \begin{equation}\label{eq: projregsol}
% \bx_d=\bV_d\underbrace{\bF_{\alpha_d}^\dagger\bff_d}_{=:\bt_d},
% \end{equation}
% where $\bV_d$ is the basis for the $d$th approximation subspace, $\bF_{\alpha_d}^{\dagger}$ is the $d$th regularized inverse for the projected problem (depending on the $d$th Tikhonov regularization parameter $\alpha_d$), and $\bff_d$ is the $d$th projected right-hand side vector. All quantities appearing in \eqref{eq: projregsol} are specific for each projection method. 
% % Depending on the projection method,  the values of $\calG(\bx_{\rm reg})$ may be different if computed for the solution $\bt_d$ of the projected problem (as defined in \eqref{eq: projregsol}) or for $\bx_d$. 
%%%
For hybrid methods based on standard Krylov subspaces, the numerator of $\calG(\bx_{\rm reg})$ (i.e., the square of the functional $\calD(\bx_{\rm reg})$) is the same when computed for both the full-dimensional and the projected problems; see \eqref{eq: hybridDP}. This is not true for GKS-based methods; see \eqref{eq: GKSDP}. Concerning the denominator of $\calG(\bx_{\rm reg})$, using the properties of the trace, one can derive the expressions
\[
\zeta - \text{trace}(\bE_d\bF_{\alpha}^\dagger),\quad\mbox{where}\quad
\bE_d=\bH_d,\bB_d,\bR_{\bA}
\]
for hybrid GMRES, hybrid LSQR and (MM)GKS, respectively. The constant $\zeta$ is $m$ for all the methods when computed for the full-dimensional problem. For the projected problems, $\zeta=d+1$ for hybrid GMRES and hybrid LSQR (see \cite{NoRu14} for more detailed derivations in the hybrid GMRES case), and $\zeta=d$ for (MM)GKS.  By default, and in agreement with the common choices made in the literature, TRIPs-Py uses the GCV criterion computed for the full-dimensional problem for hybrid GMRES and hybrid LSQR, and the projected GCV version for the solvers based on generalized Krylov subspaces; see \cite{NoRu14, bucciniGCV}. 

\subsection{Regularization operators}\label{sect:RegOp}

This section describes the regularization matrices implemented in TRIPs-Py. We consider two types of operators: those based on a finite-difference discretization of the first derivative operator and those based on framelet operators. Some more details about the usage of these regularization matrices are discussed in Section \ref{sec: TestPbs}, and an illustration is provided towards the end of the demo \verb|demo_Tomo_large_scale.ipynb|.

\paragraph{Case 1: Regularization operators based on the first derivative operator}
Let 
\begin{equation}\label{D2}
	\Psi_{D}= 
	\begin{bmatrix} 
	1 &-1 & & &  \\
	&1 &-1 & &  \\
	& &\ddots  &\ddots  & \\
	 &  & &1  &-1 
	\end{bmatrix}\in \R^{(n_D-1)\times n_D} %\in\mathbb{R}.
	\quad
	\mbox{and}\quad
	\bI_{n_D}\in\mathbb{R}^{n_D\times n_D}
\end{equation}
be a rescaled finite-difference discretization of the first derivative operator and the identity matrix of order $n_D$, respectively. For problems that depend on one or two spatial dimensions $(x,y)$ and, possibly, a time dimension $t$, the matrix $\Psi_D$ is used to obtain discretizations of the first derivatives in the $D$-direction, with $D = x$ (vertical direction), $D=y$ (horizontal direction), and $D=t$ (time direction). 
% {For simplicity, throughtout this paper, we let  $\gamma_d = 1$, but different values can be used in practice.} 
For a static image represented as a 2D array $\bX \in \R^{n_x \times n_y}$, such that $\bx = \text{vec}\left(\bX\right) \in \R^{n_x n_y}$ is obtained by stacking the columns of $\bX$, its horizontal and vertical derivatives are given as  
\begin{equation}\label{eq: Deriv}
\begin{array}{lcl}
     {\rm vec}(\Psi_{x}\bX) &=& (\bI_{n_y} \otimes \Psi_{x}) \bx \in \mathbb{R}^{n_y(n_x-1)}\\
    {\rm vec}(\bX\Psi^T_y) &=&(\Psi_{y} \otimes \bI_{n_x}) \bx \in \mathbb{R}^{(n_y -1) n_x}
\end{array},
\end{equation}
respectively. The discrete gradient is then expressed as
\begin{equation}\label{eq: gradient}
\Psi_s = [(\bI_{n_y} \otimes \Psi_{x})^T , (\Psi_{y} \otimes \bI_{n_x})^T]^T.
\end{equation}
When modelling dynamic inverse problems with a time-varying solution, let $\bX\in\R^{n_s\times n_t}$ be the 2D array whose columns store the quantity of interest at the $n_t$ time instants; note that, if such quantities are 2D images, then the columns of $\bX$ are vectorialized images with $n_s=n_xn_y$ pixels. The derivative in the time dimension is then given by 
${\rm vec}(\bX\Psi_{t}^T) =(\Psi_{t}\otimes \bI_{n_s} ) \bx \in \mathbb{R}^{(n_t -1)n_s}$. 

\paragraph{Case 2: Regularization operators based on a two-level framelet analysis operator}
When regularizing by leveraging sparsity but the desired solution is not sparse in the original domain, a well known technique is to perform a transformation to another domain, where the solution may admit a sparse representation. For this purpose, TRIPs-Py provides a two-level framelet analysis operator, 
% since it is well-known that images have sparse representation in the framelet domain. We 
defined as follows.
% \begin{dfn}\label{def: Tightframe}
Let $\bW\in\R^{r\times n}$ with $1\leq n\leq r$. The set of the rows of $\bW$ is a framelet system
for $\R^n$ if, for all $\bx\in\R^n$,
\begin{equation} \label{eq: tightframe}
\|\bx\|_2^2=\sum_{j=1}^{r}{(\bw^T_j\bx)^2},
\end{equation}
where $\bw_j\in\R^n$ denotes the $j$th row of the matrix $\bW$ (written as a column vector), 
i.e.,\linebreak[4]$\bW=[\bw_1,\bw_2,\ldots,\bw_r]^T$. The matrix $\bW$ is referred to as an analysis operator and $\bW^T$ as a synthesis operator. 
% \end{dfn}
We use the same tight frames as in \cite{buccini2020modulus, buccini2020linearized, COS09b}, i.e., the system of linear 
B-splines. This system is formed by a low-pass filter $\bW_{0}\in\mathbb{R}^{n \times n}$ 
and two high-pass filters $\bW_1,\bW_2\in\mathbb{R}^{n \times n}$, whose corresponding masks 
are
\begin{equation*}
\bw^{(0)}=\frac{1}{4}[1,2,1], \quad \bw^{(1)}=\frac{\sqrt{2}}{4}[1,0,-1], 
\quad \bw^{(2)}=\frac{1}{4}[-1,2,-1].
\end{equation*}
The analysis operator $\bW$ in one space-dimension is derived from these masks and by 
imposing reflexive boundary conditions to ensure that $\bW^T\bW=\bI$. 
The corresponding two-dimensional operator $\bW$ is given by 
\begin{equation}\label{eq: Tightframe_W}
\bW=\begin{bmatrix}
[\bW_{0}\otimes \bW_{0}]^T, 
[\bW_{0}\otimes \bW_{1}]^T,
[\bW_{0}\otimes \bW_{2}]^T,
[\bW_{1}\otimes \bW_{0}]^T,
\hdots,
[\bW_{2}\otimes \bW_{2}]^T
\end{bmatrix}^T,
\end{equation}
where $\otimes$ denotes the Kronecker product. This matrix is not explicitly formed. We
note that the evaluation of matrix-vector products with $\bW$ and $\bW^T$ is inexpensive, 
because the matrix $\bW$ is sparse. The operator $\bW$ can be used for instance as a regularization operator in GKS and MMGKS. 

% While it is important to mention that the matrices defined so far do not impose any kind of boundary conditions; nevertheless, if the behavior of the solution is known at the boundary of its domain, then one can try to recover a more accurate reconstructed approximate solution by incorporating the actual features of the boundary conditions in the regularization matrix.

\section{Overview of the TRIPs-Py test problems}\label{sec: TestPbs}
In this section, we consider three main classes of test problems. In the first class we consider both 1D and 2D deblurring, with the latter being used to produce both small-scale and large-scale synthetic test problems. {In the second class we consider computed tomography with synthetic data, which can be used to generate both small-scale and large-scale inverse problems. 
The third class contains dynamic inverse problems with real data. The usage of the classes to generate test problems and of the TRIPs-PY's solvers that can be attempted to compute their solution is illustrated in a number of demos collected in jupyter notebooks. A complete list of demos in TRIPs-Py, and short descriptions thereof, is given in Table \ref{table:demos}. 
For a smooth usage of the package, we recommend the users to first download the data needed for test problems from the google drive \footnote{\url{https://drive.google.com/drive/folders/1VB8LaFewgwNKXq7QSVns4D6Rod8W8Pa1?usp=sharing}} and place the folder data inside the folder demos.
\begin{table}
\renewcommand{\arraystretch}{1.5}
\caption{\label{table:demos} List of TRIPs-Py demos.}
\begin{tabular}{|p{64.5mm}|p{88mm}|}\hline
\textbf{Demo} & \textbf{Description}  \\ \hline\hline
\verb|demo_dynamic_Emoji.ipynb| & Emoji test problem and illustrations for solving the static and dynamic problems. \\\hline
\verb|demo_dynamic_CrossPhantom.ipynb| & CrossPhantom test problem and illustrations for solving the static and dynamic problems. \\\hline
\verb|demo_dynamic_Stempo.ipynb| & STEMPO test problem and illustrations for solving the static and dynamic problems. \\
\hline
\verb|demo_1D_Deblurring.ipynb| & 1D Deblurring test problem and illustration of direct and iterative methods to solve it.\\
\hline
\verb|demo_2D_Deblurring_small_scale.ipynb| & 2D Deblurring for a small scale test problem whose naive solution can be computed, and illustration of direct and iterative regularization methods.
\\
\hline
\verb|demo_2D_Deblurring_large_scale.ipynb|& 2D Deblurring for a large scale test problem and illustration of iterative methods to solve the ill-posed problem.
\\
\hline
\verb|demo_2D_Deblurring_your_data.ipynb| & Illustration of functions that can be used to upload users' images and define a deblurring problem.
\\
\hline
\verb|demo_Tomo_small_scale.ipynb| & Small-scale Tomography test problem and illustration of solution methods.
\\
\hline
\verb|demo_Tomo_small_scale.ipynb| & Larg-scale Tomography test problem and illustration of iterative regularization methods.
\\
\hline
\verb|demo_Tomo_saved_data.ipynb| & Illustration of how to use functions on the Tomography class that access previously saved tomography data.
\\
\hline
% \verb|demo_small_scale_CGLS.ipynb| & Illustration of CGLS on a small-scale Deblurring problem.
% \\\hline
\end{tabular}
\end{table}
\subsection{Deblurring (deconvolution)}
Deblurring can be formulated as an integral equation of the kind
\begin{equation}\label{eq: deblurring_integral_equation}
\int\mathcal{B}(s, t)x(t)ds + e = b(s),
\end{equation}
where $s, t\in \R^D$ represent spatial information (in TRIPs-Py, $D=1, 2$). The kernel $\mathcal{B}(s, t)$ (also known as the point spread function (PSF)) defines the blur. It is well known that, if the kernel is spatially invariant (as it is in TRIPs-Py), i.e., $\mathcal{B}(s, t) = \mathcal{B}(s-t)$, then \eqref{eq: deblurring_integral_equation} is a deconvolution problem. In practical settings,  discrete data are collected in finite regions, so that the continuous model \eqref{eq: deblurring_integral_equation} is  discretized and yields a linear system of equation as in \eqref{eq: linearEq}. In this particular case, the matrix $\bA\in\R^{n\times n}$ represents the blurring operator, which is defined starting from the PSF and the boundary conditions on the unknown quantity of interest. In TRIPs-Py the PSF is given by a rescaled, possibly asymmetric Gaussian of the form
\begin{equation}\label{eq: GaussPSF}
\mathcal{B}(s,t)=c\exp\left(-\frac{1}{2}(s-t)^T\bB(s-t)\right),\quad\mbox{where}\quad c>0,\;\bB=\text{diag}(\beta_1^2,\dots,\beta_D^2),
\end{equation}
whose spread parameters $\beta_1,\dots,\beta_D$ ($D=1,2$) are set by the user; reflective boundary conditions are used.  The vector $\bb$ contains the vectorized measured blurred and noisy quantity of interest. By default, the synthetic deblurring test problems available within TRIPs-Py avoid inverse crime by allowing a mismatch between the forward operator used to solve the problem, and the forward operator used to generate the data; see \cite{mueller2012linear}. Namely, the latter employs zero boundary conditions. Inverse crimes can be allowed by setting the \texttt{CommitCrime} option to \texttt{True} (the default being \texttt{False}). More details on deblurring can be found in \cite{hansen2006deblurring}. 

\subsubsection{1D Deblurring}
Referring to the notations in \eqref{eq: min}, in the 1D setting we consider a one dimensional true signal $\bx_{\rm true} \in \R^{n}$ and a convolution forward operator $\bA \in \R^{n \times n}$, with associated smoothed and noisy signal $\bb \in \R^{n}$. This problem can be setup in TRIPs-Py by defining an object of the \texttt{Deblurring1D()} class. The method $\verb|gen_xtrue()|$ generates the true signal. The required arguments are the dimension of the problem \texttt{nx} and the \texttt{test} signal; for the latter, the user can choose among saved options $\verb|piecewise|$, $\verb|sigma|$, $\verb|curve0|$, $\verb|curve1|$, $\verb|curve2|$, and $\verb|curve3|$. The user can have acces to the forward operator by calling the method $\verb|forward_Op_1D()|$ with \texttt{parameter} storing the spread of the 1D Gaussian blurring function. The data can be generated using the method $\verb|gen_data()|$. Noise is then added to the data through the method $\verb|add_noise()|$, which takes as input the noiseless data, the distribution of the random noise and the noise level (for Gaussian and Laplace noise), i.e., the ratio $\|\be\|_2/\|\bA\bx_{\rm true}\|_2$. An illustration of the usage of the class \texttt{Deblurring1D} is shown in Code 1, and more illustrations are presented in demo $\verb|demo_1D_Deblurring.ipynb|$. Figure \ref{Fig: 1DDeblurring} shows the true $\verb|'curve2'|$ signal, its blurred and noisy version with 1\% Gaussian noise, and reconstructions with tSVD and tGSVD. 
\begin{tcolorbox}[
    enhanced, attach boxed title to top center={yshift=-2mm},
    colback=cyan!10,
    colframe=black,
    colbacktitle=black,
    title = Code 1: Generate Image deblurring 1D,
    text width = 15cm,
    fonttitle=\bfseries\color{white},
    boxed title style={size=small,colframe=darkspringgreen,sharp corners},
    sharp corners,]
% \begin{minted}[breaklines]{python}
\begin{lstlisting}[language=Python]
# Define an object of Deblurring1D class
Deblur1D = Deblurring1D(CommitCrime = True)
nx = 200
# Generate a true solution
x_true = Deblur1D.gen_xtrue(nx, test = 'curve2')
# Generate the forward operator
A = Deblur1D.forward_Op_1D(parameter = 30, nx = nx) 
# Generate the blurred and noisy curve
b_true = Deblur1D.gen_data(x_true)
(d, delta) = Deblur1D.add_noise(b_true, 'Gaussian', noise_level = 0.01)
Deblur1D.plot_data(d)
Deblur1D.plot_rec(x_true)
\end{lstlisting}
% \end{minted}
\end{tcolorbox}

% \begin{tcolorbox}[
%     enhanced, attach boxed title to top center={yshift=-2mm},
%     colback=cyan!15,
%     colframe=black,
%     colbacktitle=black,
%     title= Code 1: Generate Image deblurring 1D,
%     text width = 15cm,
%     fonttitle=\bseries\color{white},
%     boxed title style={size=small,colframe=darkspringgreen,sharp corners},   sharp corners,]
% \begin{minted}{python}
% # Define an object of Deblurring1D class
% Deblur1D = Deblurring1D(CommitCrime = True)
% nx = 200
% # Generate a true solution
% x_true = Deblur1D.gen_xtrue(nx, test = 'curve2')
% # Generate the forward operator
% A = Deblur1D.forward_Op_1D(parameter = 30, nx = nx) 
% # Generate the blurred and noisy curve
% b_true = Deblur1D.gen_data(x_true)
% (d, delta) = Deblur1D.add_noise(b_true, 'Gaussian', noise_level = 0.01)
% Deblur1D.plot_data(d)
% Deblur1D.plot_rec(x_true)
% \end{minted}
% \end{tcolorbox}
\begin{figure}[h!]
\centering
\begin{tabular}{cccc}
	% \begin{minipage}{0.19\textwidth}
	\includegraphics[width=0.22\textwidth]{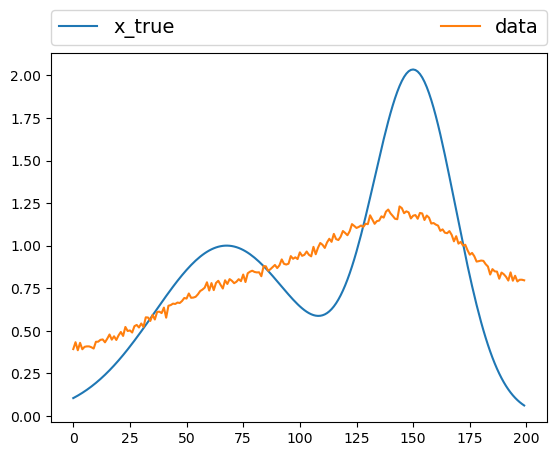}&
	% \end{minipage}
	% \begin{minipage}{0.19\textwidth}
	\includegraphics[width=0.22\textwidth]{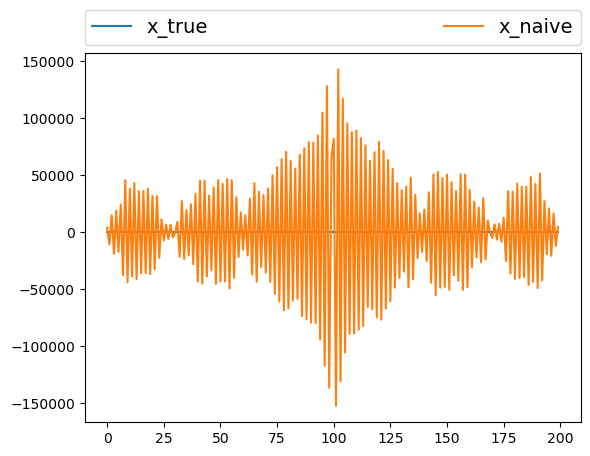}&
	% \end{minipage}
	% \begin{minipage}{0.19\textwidth}
	% 	\includegraphics[width=0.17\textwidth]{Figures/1DDeblurring/1D_TSVD.png}&
	% % \end{minipage}
 % % \begin{minipage}{0.19\textwidth}
		\includegraphics[width=0.22\textwidth]{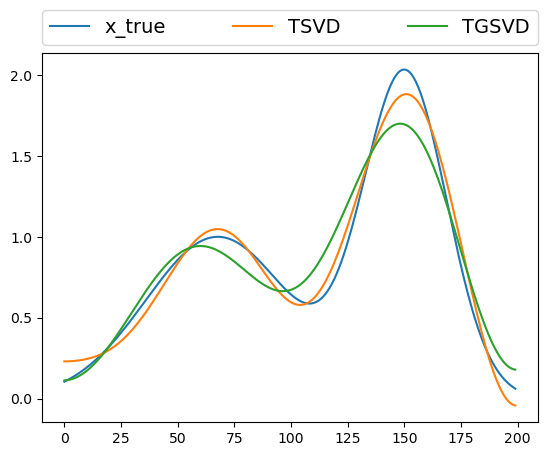}&
	% \end{minipage}
  % \begin{minipage}{0.19\textwidth}
		\includegraphics[width=0.22\textwidth]{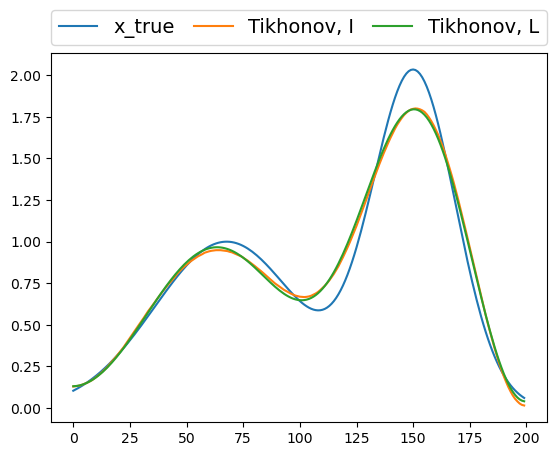}\\
	% \end{minipage}\\
(a) & (b)  & (c)  & (d) 
\end{tabular}
	\caption{1D image deblurring test problem. (a) $\bx_{\rm true}$ and data with $1\%$ Gaussian noise, (b) $\bx_{\rm true}$ and naive, unregularized solution, (c) $\bx_{\rm true}$, TSVD and TGSVD solutions, (d) $\bx_{\rm true}$, standard and general form Tikhonov solutions.}
	\label{Fig: 1DDeblurring}
\end{figure}
\subsubsection{2D deblurring}

In this section we show strategies to define and solve a 2D deblurring problem in TRIPs-Py. Similarly to the 1D case, 2D image deblurring test problems are set up as objects of the \texttt{Deblurring2D()} class. As a first step, we start with a small scale problem to illustrate the ill-posedness of such inverse problem and the performance of direct regularization methods applied to it (with some more illustrations available through the notebook $\verb|demo_Deblurring_small_scale.ipynb|$). Then we consider a larger scale problem to illustrate iterative regularization methods (a full investigation is provided in the demo $\verb|demo_Deblurring_large_scale.ipynb|$). For consistancy, we choose the same image in both examples, but the user can choose from other images provided in the package such as $\verb|h_im|$, $\verb|hubble|$, $\verb|grain|$, and $\verb|sky|$. Other images (in tiff, tif, jpg, png or mat formats) can be easily used in TRIPs-Py by performing the following operations:
    \begin{enumerate}
        \item Create the folder \verb|my_image_data| under demos/data and place the desired image inside the folder.
        \item Run the function\vspace{-0.2cm}

\begin{tcolorbox}[
    enhanced, attach boxed title to top center={yshift=-2mm},
    colback= cyan!3,
    colframe=white,
    colbacktitle=black,
    % title= Code 2: Generate Image deblurring 2D,
    text width = 15cm,
    fonttitle=\bseries\color{white},
    boxed title style={size=small,colframe=darkspringgreen,sharp corners},
    sharp corners,]
% \begin{minted}[breaklines]{python}
\begin{lstlisting}[language=Python]
convert_image_for_trips(imag = 'image_name', image_type= image_type')
\end{lstlisting}
\end{tcolorbox}
        % \[\verb|convert_image_for_trips(imag = `image_name', image_type= `image_type')|,\vspace{-0.2cm}\] 
where $\verb|'image_name'|$ and $\verb|'image_type'|$ refer to the desired image name and type.
    \end{enumerate}
This is shown in  the demo $\verb|demo_Deblurring_your_data.ipynb|$.

% , described below. In TRIPs-Py there are three demos that illustrate the 2D image deblurring problem as below.

% \begin{enumerate}
%     \item [\diamond] $\verb|demo_Deblurring_small_scale.ipynb|$
%     A jupyter notebook that defines a relatively small scale deblurring problem and illustrates regularization methods designed for small scale problems such as TSVD and Tikhonov.
%     \item [\diamond] $\verb|demo_Deblurring_large_scale.ipynb|$ In this notebook we illustrate how to define a large scale deblurring problem along with illustration of iterative regularization methods such as hybrid based methods (Hybrid\_GMRES and Hybrid\_LSQR) and MMGKS.
%     \item [\diamond]  $\verb|demo_Deblurring_your_data.ipynb|$
%     On this notebook we illustrate how the user can handle their own data. 
%   When the data in one of the following formats (tiff, tif, jpg, png or mat) the user can follow the instructions probided below to use them in TRIPs-Py.
%     \begin{enumerate}
%         \item Step 1: Create the folder my\_image\_data under demos/data and place your image.
%         \item Run the function convert\_image\_for\_trips(imag = `image\_name', image\_type= `image\_type')
%     \end{enumerate}
% After running steps (a) and (b) above, a file `image\_name.mat' will be stored on your permanent image folder image\_data which is in compatible for to be used from functionalities in TRIPs-Py.
% \end{enumerate}
\paragraph{Small-scale image deblurring} For this illustration we set $n_x= n_y = 50$ and define the forward operator to be a Gaussian PSF \eqref{eq: GaussPSF} with parameters $(\beta_1,\beta_2)=(1,1)$. The true image and the blurred and noisy image with $1\%$ Gaussian noise are shown in Figure \ref{Fig: smallScaleDeblurring} a) and b). We call methods TSVD, Hybrid LSQR, and MMGKS with parameters specified as follows. 
% \[
% \verb|(x_tsvd, truncation_value) = TruncatedSVD_sol(A.todense(), b_vec, regparam = 'dp', delta = delta)|
% \]
% where \texttt{delta} is the noise level (returned from \texttt{add\_noise}). 

\begin{tcolorbox}[
    enhanced, attach boxed title to top center={yshift=-2mm},
    colback= cyan!3,
    colframe=white,
    colbacktitle=black,
    % title= Code 2: Generate Image deblurring 2D,
    text width = 15cm,
    fonttitle=\bseries\color{white},
    boxed title style={size=small,colframe=darkspringgreen,sharp corners},
    sharp corners,]
% \begin{minted}[breaklines]{python}
\begin{lstlisting}[language=Python]
(x_tsvd, truncation_value) = TruncatedSVD_sol(A.todense(), b_vec, regparam = 'dp', delta = delta)
(x_hybrid_lsqr, info_hybrid_lsqr) = Hybrid_LSQR(A, b_vec, n_iter = 100, regparam = 'dp', x_true = x_true, delta = delta)
L = first_derivative_operator_2d(nx, ny)
(x_mmgks, info_mmgks) = MMGKS(A, data_vec, L, pnorm=2, qnorm=1, projection_dim=2, n_iter = 100, regparam = 'dp', x_true = x_true, delta = delta)
\end{lstlisting}
\end{tcolorbox}
We remark that the user should provide the operator $\bA$ as a dense matrix for TSVD and the observed data as a vector for all the methods. The regularization parameter can be a scalar, or chosen by the discrepancy principle (\verb|'dp'|) or generalized cross validation (\verb|'gcv'|). When the regularization parameter is set to \verb|'dp'|, the user must specify the noise level (returned from \texttt{add\_noise()} function (see Code 1)). For the MMGKS, the user can specify the regularization parameter $\bL$ (for this example we choose $\bL$  to be 2D discretized derivative operator of the first order). The values of $p$ and $q$ can be set as well through the input parameters \verb|pnorm| and \verb|qnorm|. 
Approximate solutions obtained from truncated SVD, and iterative methods such as Hybrid LSQR and MMGKS are shown in Figure \ref{Fig: smallScaleDeblurring} c), d), and e), respectively. The approximate solution and information collected through the iterations is outputed.
Apart from visual inspections provided in Figure \ref{Fig: smallScaleDeblurring} for this example, more quantitative measures on the reconstructed solution can be found on the demo jupyter notebook $\verb|demo_Deblurring_small_scale.ipynb|$.

\begin{figure}[h!]
\centering
	\centering
\begin{tabular}{ccccc}
	% \begin{minipage}{0.19\textwidth}
		\includegraphics[width=0.17\textwidth]{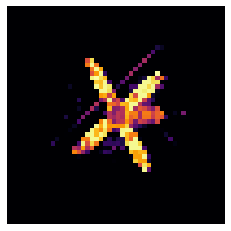}&
	% \end{minipage}
	% \begin{minipage}{0.19\textwidth}
		\includegraphics[width=0.17\textwidth]{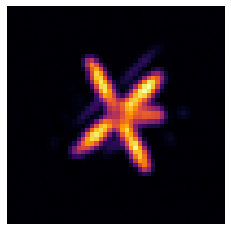}&
	% \end{minipage}
	% \begin{minipage}{0.19\textwidth}
		\includegraphics[width=0.17\textwidth]{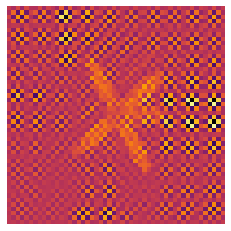}&
	% \end{minipage}
 % \begin{minipage}{0.19\textwidth}
		\includegraphics[width=0.17\textwidth]{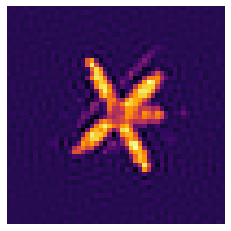}&
	% \end{minipage}
  % \begin{minipage}{0.19\textwidth}
		\includegraphics[width=0.17\textwidth]{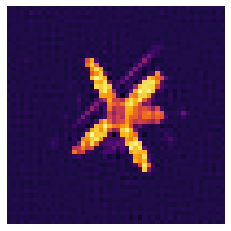}\\
	% \end{minipage}\\
(a) & (b)  & (c)  & (d) & e)
\end{tabular}
	\caption{2D image deblurring small-scale test problem. (a) True image of $50 \times 50$ pixels, (b) Blurred and noisy image with $1\%$ Gaussian noise. Approximate solutions by (c) TSVD, (d) Hybrid LSQR, and (e) MMGKS.}
	\label{Fig: smallScaleDeblurring}
\end{figure}

\paragraph{Large-scale image deblurring}
For the large-scale version we consider a satellite image of size $128\times 128$ pixels shown in Figure \ref{Fig: largeScaleDeblurring} (a). Such image is blurred by a Gaussian PSF with parameters $(\beta_1,\beta_2)=(3,3)$ and we add $1\%$ Gaussian noise. An illustration of how to define a 2D Deblurring problem is shown in Code 2. The blurred and noisy image is shown in Figure \ref{Fig: largeScaleDeblurring} (b). Approximate solutions obtained by hybrid GMRES, hybrid LSQR, and MMGKS are shown in Figure \ref{Fig: largeScaleDeblurring} c), d), and e). The calls to these solvers are essentially identical to the ones illustrated for the small-scale example. More details can be found on demo $\verb|demo_Deblurring_large_scale.ipynb|$. 
\begin{figure}[h!]
\centering
	\centering
\begin{tabular}{ccccc}
	% \begin{minipage}{0.19\textwidth}
		\includegraphics[width=0.17\textwidth]{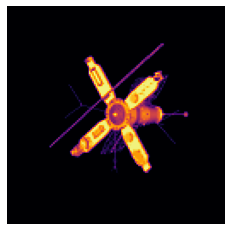}&
	% \end{minipage}
	% \begin{minipage}{0.19\textwidth}
		\includegraphics[width=0.17\textwidth]{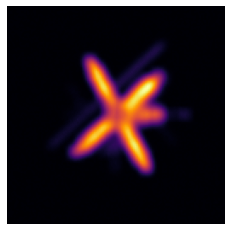}&
	% \end{minipage}
	% \begin{minipage}{0.19\textwidth}
		\includegraphics[width=0.17\textwidth]{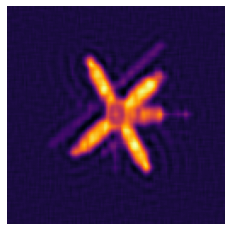}&
	% \end{minipage}
 % \begin{minipage}{0.19\textwidth}
		\includegraphics[width=0.17\textwidth]{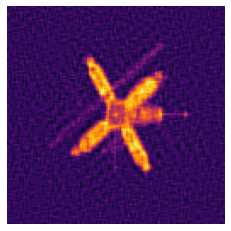}&
	% \end{minipage}
  % \begin{minipage}{0.19\textwidth}
		\includegraphics[width=0.17\textwidth]{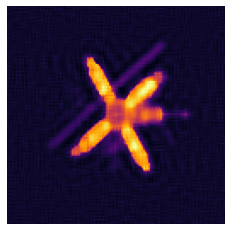}\\
	% \end{minipage}\\
(a) & (b)  & (c)  & (d) & e)
\end{tabular}
	\caption{2D image deblurring large-scale test problem. (a) True image of $128 \times 128$ pixels. (b) Blurred and noisy image with $1\%$ Gaussian noise. Approximate solutions obtained by (c) Hybrid GMRES, (d) Hybrid LSQR, and (e) MMGKS.}
	\label{Fig: largeScaleDeblurring}
\end{figure}

  \begin{tcolorbox}[
    enhanced, attach boxed title to top center={yshift=-2mm},
    colback=cyan!10,
    colframe=black,
    colbacktitle=black,
    title = Code 2: Generate Image deblurring 2D,
    text width = 15cm,
    fonttitle=\bfseries\color{white},
    boxed title style={size=small,colframe=darkspringgreen,sharp corners},
    sharp corners,]
\begin{lstlisting}[language=Python]
% \begin{minted}[breaklines]{python}
# Create an object of class Deblurring()
Deblur = Deblurring(CommitCrime = False)
nx = 256
ny = 256
# Generate the forward operator
A = Deblur.forward_Op((9,9), (3,3), nx, ny)
choose_image = 'satellite'
# generate the true image
x_true = Deblur.gentrue(choose_image)
# generate data
b_true = Deblur.gen_data(x_true) # generate_matrix)
# add noise
(b, delta) = Deblur.add_noise(b_true, opt= 'Gaussian', noise_level = 0.01)
b_vec = b.reshape((-1,1))
\end{lstlisting}
\end{tcolorbox}

% \begin{tcolorbox}[
%     enhanced, attach boxed title to top center={yshift=-2mm},
%     colback=cyan!15,
%     colframe=black,
%     colbacktitle=black,
%     title= Code 2: Generate Image deblurring 2D,
%     text width = 15cm,
%     fonttitle=\bseries\color{white},
%     boxed title style={size=small,colframe=darkspringgreen,sharp corners},
%     sharp corners,]
% \begin{minted}{python}
% # Create an object of class Deblurring()
% Deblur = Deblurring(CommitCrime = False)
% nx = 256
% ny = 256
% # Generate the forward operator
% A = Deblur.forward_Op((9,9), (3,3), nx, ny)
% choose_image = 'satellite'
% # generate the true image
% x_true = Deblur.gentrue(choose_image)
% # generate data
% b_true = Deblur.gen_data(x_true) # generate_matrix)
% # add noise
% (b, delta) = Deblur.add_noise(b_true, opt= 'Gaussian', noise_level = 0.01)
% b_vec = b.reshape((-1,1))
% \end{minted}
% \end{tcolorbox}

\subsection{Computerized tomography}
Another test problem considered in TRIPs-Py is 2D X-ray computerized tomography (CT), which consists in reconstructing an object (i.e., the attenuation coefficients of an object) from a set of projections along straight lines (i.e., intensities of energy rays recorded by detectors). In the following we briefly describe the physical and mathematical formulation of CT and then illustrate how to define a CT test problem within TRIPs-Py. 

Consider $\bx=[x_1,x_2]^T$, $\beta \in [-\pi, \pi]$, $\boldsymbol{\beta}=[\cos(\beta), \sin(\beta)]^T$, $\boldsymbol{\beta}^{\perp}=[-\sin(\beta), \cos(\beta)]^T$, and $c\in \R$. Let
\[
\mathcal{L}(\beta, c) =\left\{\bx \in \R^2\,|\,\bx\cdot\boldsymbol{\beta} = c \right\}
=\left\{\bx \in \R^2\,|\,\bx=c\boldsymbol\beta+\ell\boldsymbol\beta^\perp,\:\ell\in\R\right\}
\]
be the perpendicular line to $\boldsymbol\beta$ with signed orthogonal distance $c$ from the origin. Then, assuming that absorption dominates potential scattering effects, Lambert-Beer law links the attenuation coefficient $f(\bx)$ of the object we wish to image to the recorded intensity $I_{\beta,c}$ of the measured X-ray of incoming intensity $I_0$ along $\mathcal{L}(\beta, c)$ as follows
\begin{equation}\label{eq: CT}
\underbrace{\int_{\mathcal{L}(\beta, c)}f(\bx)d\ell}_{=\mathfrak{R}[f](\beta, c)}=
% \underbrace{-\log\left(\frac{I_{\beta,s}}{I_0}+e_{\beta,c}\right)}_{=b(\beta,c)}.
-\log\left(\frac{I_{\beta,s}}{I_0}+e_{\beta,c}\right),
\end{equation}
where $e_{\beta,c}$ is a random perturbation corrupting the measurements. The left-hand side of the above equation, when computed for all $\beta$'s and $c$'s, is the Radon transform of the function $f$.
% We let $f$ to be the density of the medium whose interior is of interest. 
% Beer's Law \cite{BL} gives the data measured from an X-ray CT scan as $\mathfrak{R}f(\theta, c) =\int_{\bx \in \mathcal{L}(\theta, c)} f(\bx)d\bs$ with $ds$ representing the arc length measure. 
After a discretization process, from \eqref{eq: CT} we obtain a discrete formulation of the problem, $\bA\bx = \bb$. More details on CT can be found in \cite{PCHCT}.

Similarly to deblurring, we generate a CT test problem in TRIPs-Py by first defining an object of the class \texttt{Tomography()}. 
By default, also the synthetic CT test problems available within TRIPs-Py avoid inverse crime by using slightly different forward operators to solve the problem and to generate the data. Namely, the set of projection angles for the two operators are affected by a small constant mismatch. Inverse crimes can be allowed by setting the \texttt{CommitCrime} option to \texttt{True} (the default being \texttt{False}).
The following demos in TRIPs-Py illustrate how to set up and solve tomography test problems.
\begin{enumerate}
    \item [$\diamond$] $\verb|demo_Tomo_small_scale.ipynb|$
    This notebook defines a tomography test problem with small dimensions, so that the forward operator can be explicitly formed and stored to compute a naive solution or a regulatized solution from truncated SVD. 
    \item [$\diamond$] $\verb|demo_Tomo_large_scale.ipynb|$ This notebook showcases how to generate a large-scale tomography problem along with how to call algorithms to solve the formulated inverse problem. 
    For both demos described above, to setup the forward operator, we use the ASTRA toolbox with `fanflat' 2D geometry. Other geometries such as `parallel' and `fanflat\_vec' can be used. More details on the ASTRA toolbox and its documentation can be found in \cite{van2015astra} and references therein.
    \item [$\diamond$]  $\verb|demo_Tomo_saved_data.ipynb|$
    This notebook demonstrates how to generate a tomogrpahy problem where the forward operator and the data are already available. This demo has no dependency on any ASTRA toolbox to generate the forward operator. 
\end{enumerate}
% MENTION PARALLEL LINE GEOMETRY AND DISCRETIZATION PROCESS, AND ASTRA.

Code 3 illustrates how to set up a large scale tomography test problem, where $n_x = n_y = 256$ and the number of view angles is 50; a small scale problem can be defined similarly by reducing $n_x$ and $n_y$. The angles can be limited by varying the parameter $\verb|views|$. The true phantom is given in Figure \ref{Fig: Tomography}(a) and the observed sinogram with $1\%$ Gaussian noise is shown in Figure \ref{Fig: Tomography}(b). Reconstructed solutions by iterative methods Hybrid LSQR, GKS, and MMGKS are shown in Figure \ref{Fig: Tomography} (c), (d), and (e), respectively.

  \begin{tcolorbox}[
    enhanced, attach boxed title to top center={yshift=-2mm},
    colback=cyan!10,
    colframe=black,
    colbacktitle=black,
    title = Code 3: Generate X-Ray CT problem,
    text width = 15cm,
    fonttitle=\bfseries\color{white},
    boxed title style={size=small,colframe=darkspringgreen,sharp corners},
    sharp corners,]
% \begin{minted}[breaklines]{python}
\begin{lstlisting}[language=Python]
# Define an object of the class Tomogrpahy
Tomo = Tomography()
# Specify the dimensions of the phantom and the number of angels
nx = 256
ny = 256
views = 50
# Define the true solution throught the function gen_true(). 
(x_true, nx, ny) = Tomo.gen_true(testproblem = 'tectonic', nx = nx, ny = ny)
# Define the forward operator
A = Tomo.forward_Op(nx, ny, views)
# Generate the data
(A, b_true, p, q, AforMatrixOperation) = Tomo.gen_data(x_true, nx, ny, views)
# Add noise in the true simulated sinogram b_true
(b, delta) = Tomo.add_noise(b_true = b_true, opt = 'Gaussian', noise_level = 0.001)
\end{lstlisting}
\end{tcolorbox}

% \begin{tcolorbox}[
%     enhanced, attach boxed title to top center=r{yshift=-2mm},
%     colback=cyan!15,
%     colframe=black,
%     colbacktitle=black,
%     title = Code 3: Generate X-Ray CT problem,
%     text width = 15cm,
% fonttitle=\bseries\color{white},
%     boxed title style={size=small,colframe=darkspringgreen,sharp corners},
%     sharp corners,]
% \begin{minted}{python}
% # Define an object of the class Tomogrpahy
% Tomo = Tomography()
% # Specify the dimensions of the phantom and the number of angels
% nx = 256
% ny = 256
% views = 50
% # Define the true solution throught the function gen_true(). 
% (x_true, nx, ny) = Tomo.gen_true(testproblem = 'tectonic', nx = nx, ny = ny)
% # Define the forward operator
% A = Tomo.forward_Op(nx, ny, views)
% # Generate the data
% (A, b_true, p, q, AforMatrixOperation) = Tomo.gen_data(x_true, nx, ny, views)
% # Add noise in the true simulated sinogram b_true
% (b, delta) = Tomo.add_noise(b_true = b_true, opt = 'Gaussian', noise_level = 0.001)
% \end{minted}
% \end{tcolorbox}

\begin{figure}[h!]
\centering
	\centering
\begin{tabular}{ccccc}
		\includegraphics[width=0.17\textwidth]{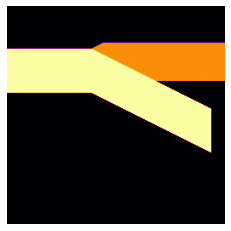}&
		\includegraphics[width=0.17\textwidth]{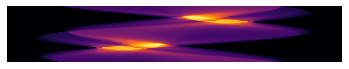}&
		\includegraphics[width=0.17\textwidth]{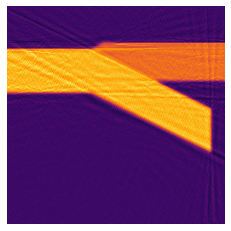}&
		\includegraphics[width=0.17\textwidth]{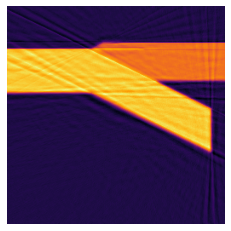}&
		\includegraphics[width=0.17\textwidth]{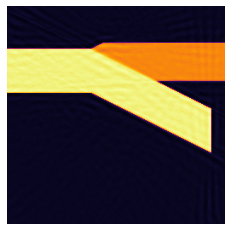}\\
(a) & (b)  & (c)  & (d) & e)
\end{tabular}
	\caption{Tomography test problem. (a) True image of $256 \times 256$ pixels. (b) Sinogram. Approximate solutions obtained by (c) Hybrid LSQR, (d) GKS, and (e) MMGKS.}
	\label{Fig: Tomography}
\end{figure}

\subsection{Dynamic Computerized Tomography}\label{sec: DemoDynamic}
Recent technological advancements in detector speed and accuracy 
%multichannel photon counting has caused a 
resulted in a growing interest in X-ray computerized tomography, which in turn prompted the need for new and efficient methods to analyze the collected data. In particular, in this paper and within TRIPs-Py, we are interested in reconstructing a series of images to explore the spatial and temporal properties of data acquired within dynamic CT.
% , to monitor structural evolution or movement during the scanning process. 
The aim of this section is three-folded: 1) to briefly describe several regularization terms $\mathcal{R}\left(\bx\right)$ for use in dynamic inverse problems, 2) to show how to handle dynamic data and generate dynamic CT problems in TRIPs-Py, and 3) to illustrate the need for temporal regularization along with spatial regularization, as well as the usage of TRIPs-Py solvers for dynamic CT.

In the dynamic setting, we are interested in solving the minimization problem

\begin{equation}\label{eq: blockF}
\underbrace{
\begin{bmatrix} \widehat\bA^{(1)} 
\\ &  \ddots\\
& & \widehat\bA^{(n_t)} 
\end{bmatrix}}_{\bA}
\underbrace{\begin{bmatrix}
    \widehat\bx^{(1)}\\
    \vdots\\
    \widehat\bx^{(n_t)}
\end{bmatrix}}_{\bx}
+
\underbrace{\begin{bmatrix}
    \widehat\be^{(1)}\\
    \vdots\\
    \widehat\be^{(n_t)}
\end{bmatrix}}_{\be}
=
\underbrace{
\begin{bmatrix}
    \widehat\bb^{(1)}\\
    \vdots\\
    \widehat\bb^{(n_t)}
\end{bmatrix}}
_{\bb}.
\end{equation}
The operator $\bA$ has typically a block diagonal structure where the blocks $\widehat\bA^{(t)}$ change in time $t=1,\dots,n_t$. If 
%the operator is the same at every time instance (i.e., 
$\widehat\bA^{(t)} = \widehat\bA$ for $t=1,\dots,n_t$, the forward operator simplifies to $\bA = \bI_{n_t} \otimes \widehat\bA$, where $\otimes$ denotes the Kronecker product. 
% Although \eqref{eq: blockF} may be approached by solving a sequence of `static' regularized problems of the form
% \begin{equation}\label{eq: static}
% \bx^{(t)}_\text{static} = \argmin_{\bx\in \R^{ns}} \| \bA^{(t)}\bx -\bb^{(t)}\|_2^2 + \alpha \calR(\bx), \qquad t = 1,2,\dots, n_t, 
% \end{equation}
% often, in order to overcome the difficulties associated to having only limited information per time instance and to improve the reconstruction quality, it is crucial to incorporate temporal information into the reconstruction process by solving 
% \begin{equation}\label{eq: l2-lq}
% \bx_\text{dynamic}=\argmin_{\bx \in \R^{n}}\|\bA\bx-\bb\|_2^2 + \alpha\calR(\bx).
% \end{equation}
Problem \eqref{eq: blockF} may be approached by solving a sequence of `static' regularized problems of the form
\begin{equation}\label{eq: static}
\bx^{(t)}_\text{static} = \argmin_{\bx\in \R^{n_s}} \| \widehat\bA^{(t)}\bx -\widehat\bb^{(t)}\|_2^2 + \alpha \|\Psi_s\bx\|_q^q, \qquad t = 1,2,\dots, n_t, 
\end{equation}
where $\Psi_s$ is defined as in \eqref{eq: gradient} and $q>0$ (note that $q=1$ corresponds to anisotropic total variation in space). However, in order to overcome the difficulties associated to often having only limited information per time instance, and to improve the reconstruction quality, it is crucial to incorporate temporal information into the reconstruction process by solving 
\begin{equation}\label{eq: l2-lq}
\bx_\text{dynamic}=\argmin_{\bx \in \R^{n}}\|\bA\bx-\bb\|_2^2 + \alpha\calR(\bx).
\end{equation}
Here $\calR(\bx)$ is a regularization term that takes into consideration both spatial and temporal dimensions.
% where $\calR(\bx)$ denotes the regularization term and in its general form is set to $\|\Psi \bx\|_q^q$, with $q > 0$ and $\Psi = \bL_s$.
% Moreover, in the dynamic inverse problems setting, another difficulty comes from the limited information available for
% time instance during the measurement process. For such case, temporal information is crucial to be incorporated into the regularization to improve reconstruction quality. 
% Hence, we seek to solve 
%
In TRIPs-Py, the user can find three regularization terms for \eqref{eq: l2-lq}; more details and a more extended list of such regularizers can be found in \cite{pasha2021efficient}. 
In the Bayesian setting, edge-preserving regularization methods for dynamic problems are proposed in \cite{lan2023spatiotemporal}. 
\begin{itemize}
\item [$\diamond$] \textbf{Anisotropic space-time TV} We take
\begin{equation}\label{eqn:TVan}
\calR(\bx)= \>  \sum_{t=1}^{n_t} \|\Psi_s\bx^{(t)}\|_1 +  \sum_{t=1}^{n_t-1}\|\bx^{(t+1)}-\bx^{(t)}\|_1\\
= \> \|(\bI_{n_t}\otimes\Psi_s)\bx\|_1 + \|(\Psi_t\otimes \bI_{n_s})\bx\|_1%\\,
\end{equation}
The anisotropic TV terms $\|\Psi_s\bx^{(t)}\|_1$,  $t=1,\dots, n_t$, ensure that the discrete spatial gradients of the images are sparse at each time step; moreover, we enforce that the images do not change considerably from one time instant to the next one by penalizing the 1-norm of their difference.
\item[$\diamond$] \textbf{Isotropic TV in space, Anisotropic TV in time}
Assuming, for simplicity, that $n_x=n_y$, we take
% \begin{equation}\label{eqn:TViso}
% \begin{aligned}
% \rm TV_{\rm iso}(\bx) &= \>  \sum_{\ell=1}^{nx \cdot ny \cdot nt}\sqrt{
% ({\bbz}_y(\bx))_{\ell}^2 + ({\bbz}_x(\bx))_{\ell}^2} +  \sum_{t=1}^{nt-1}\|\bx^{(t+1)}-\bx^{(t)}\|_1\\
% &= \> \|\,[{\bbz}_y(\bu), {\bbz}_x(\bx)]\,\|_{2,1} + \|(\Psi_t\otimes \bI_{ns})\bx\|_1,
% \end{aligned} 
% \end{equation}
% where $\|\cdot\|_{2,1}$ denotes the functional defined, for a matrix  $\bY
% \in\mathbb{R}^{my\times ny}$, as $
% \|\bY\|_{2,1}=\sum_{i=1}^{my}\;
% \|\bY_{i,:}\|_2$ and
% \begin{equation}\label{eq:z3DisoTV}
% \begin{array}{ccl}
% {\bbz}_y(\bx):=&\>(\bI_{nt}\otimes\bI_{nx}\otimes{\bar{\Psi}}_{y})\bx\,,\\
% {\bbz}_x(\bx):=&\>(\bI_{nt}\otimes{\bar{\Psi}}_{h}\otimes\bI_{nx})\bx.\\
% \end{array}
% \end{equation}
\begin{equation}\label{eqn:TViso}
\begin{aligned}
\calR(\bx)&= \>  \sum_{\ell=1}^{n}\sqrt{
((\bI_{n_t}\otimes\bI_{n_y}\otimes{\bar{\Psi}}_{x})\bx)_{\ell}^2 + ((\bI_{n_t}\otimes{\bar{\Psi}}_{y}\otimes\bI_{n_x})\bx)_{\ell}^2} +  \sum_{t=1}^{n_t-1}\|\bx^{(t+1)}-\bx^{(t)}\|_1\\
&= \> \|\,[(\bI_{n_t}\otimes\bI_{n_y}\otimes{\bar{\Psi}}_{x})\bx, (\bI_{n_t}\otimes{\bar{\Psi}}_{y}\otimes\bI_{n_x})\bx]\,\|_{2,1} + \|(\Psi_t\otimes \bI_{n_s})\bx\|_1,
\end{aligned} 
\end{equation}
where ${\bar{\Psi}}_{D}$ ($D=x, y$) denote the square version of the matrix $\Psi_D$ defined in \eqref{D2}, where a zero row has been added at the bottom, and $\|\cdot\|_{2,1}$ denotes the functional defined, for a matrix  $\bZ
\in\mathbb{R}^{m_x\times m_y}$, as $
\|\bZ\|_{2,1}=\sum_{i=1}^{m_x}\;
\|\bZ_{i,:}\|_2$. 
% and
% \begin{equation}\label{eq:z3DisoTV}
% \begin{array}{ccl}
% {\bbz}_y(\bx):=&\>(\bI_{nt}\otimes\bI_{nx}\otimes{\bar{\Psi}}_{y})\bx\,,\\
% {\bbz}_x(\bx):=&\>(\bI_{nt}\otimes{\bar{\Psi}}_{h}\otimes\bI_{nx})\bx.\\
% \end{array}
% \end{equation}
\item [$\diamond$] \textbf{Group sparsity} 
% Group sparsity allows to promote sparsity when reconstructing a vector of unknown pixels that are naturally partitioned in subsets; see~\cite{bach2012optimization}. 
In dynamic CT, one can naturally group the spatial variables (pixels) at each time instant, i.e., $\{\bx^{(t)}\}_{t=1}^{n_t}$, although there are other possible ways of defining groups. In TRIPs-Py we impose group sparsity across the groups defined by the pixels of the gradient images for all the time instants, i.e., we take 
\begin{equation}\label{eq: SparsityReg} 
\calR(\bx) = \sum_{\ell=1}^{n_s'} \left( \sum_{t=1}^{n_t} (\Psi_s\bx^{(t)})_{\ell}^2 \right)^{1/2} = \|\Psi_s\bX\|_{2,1},
\end{equation}
where $n_s' = (n_y-1)n_x + (n_y-1)n_x$ is the total number of pixels in the gradient images and, referring to the notations in \eqref{eq: blockF}, $\bX=[\widehat\bx^{(1)},\dots,\widehat\bx^{(n_t)}]\in\R^{n_s\times n_t}$. 
%
% \begin{equation}\label{eq:Zgroup}
% \begin{aligned}
% \bZ &=[\Psi_s\bx^{(1)},\dots,\Psi_s\bx^{(n_t)}]=\Psi_s\bX\in\R^{ns'\times nt},\\
% \bz &={\rm{vec}}(\bZ)=(\bI_{nt}\otimes \Psi_s)\bx\,.
% \end{aligned}
% \end{equation}
%
% To enforce piecewise constant structure in space and time, we adopt the following approach. Let $n_s' = (n_y-1)n_x + (n_y-1)n_x$ be the total number of pixels in the gradient images. 
% Consider the groups defined by the vectors
% \[ \bz_{\ell} = \begin{bmatrix} (\Psi_s\bx^{(1)})_{\ell},\dots,(\Psi_s\bx^{(nt)})_{\ell} \end{bmatrix} = \left(\bI_{nt} \otimes \be_\ell^T\Psi_s\right)\bx \in \mathbb{R}^{nt}, \qquad \ell=1,\dots,ns'.\]
% Alternatively, define the matrix $\bZ$ whose columns represent the vectorized gradient images at different time $t$ as 
% \begin{equation}\label{eq:Zgroup}
% \begin{aligned}
% \bZ &=[\Psi_s\bx^{(1)},\dots,\Psi_s\bx^{(nt)}]=\Psi_s\bX\in\R^{ns'\times nt},\\
% \bz &={\rm{vec}}(\bZ)=(\bI_{nt}\otimes \Psi_s)\bx\,.
% \end{aligned}
% \end{equation}
% The regularization term corresponding to group sparsity can then be expressed as a mixture of norms 
% \begin{equation}\label{eq: SparsityReg} 
% \rm GS(\bx) := \sum_{\ell=1}^{n_s'}\|\bz_{\ell}\|_2 = \sum_{\ell=1}^{n_s'} \left( \sum_{t=1}^{nt} (\Psi_s\bx^{(t)})_{\ell}^2 \right)^{1/2} = \|\Psi_s\bX\|_{2,1}. \nonumber 
% \end{equation}
\end{itemize}

Code 4 gives examples of the usage within TRIPs-Py of the basic regularization operators that appear in the regularization functionals listed above, as well as the framelet operator described in Section \ref{sect:RegOp}. These are coded in the file \verb|operators.py| under the directory \verb|trips/utilities|. All such regularization functionals can be handled in the general framework of MMGKS, as described at the end of Section \ref{sec: GKS}.
  \begin{tcolorbox}[
    enhanced, attach boxed title to top center={yshift=-2mm},
    colback=cyan!10,
    colframe=black,
    colbacktitle=black,
    title = Code 4: Regularization operators,
    text width = 15cm,
    fonttitle=\bfseries\color{white},
    boxed title style={size=small,colframe=darkspringgreen,sharp corners},
    sharp corners,]
% \begin{minted}[breaklines]{python}
\begin{lstlisting}[language=Python]
# Discretization of the first derivative operator in 1D
L = gen_first_derivative_operator(n)
# Discretization of the first derivative operator in 2D
L = gen_first_derivative_operator_2D(nx, ny)
# Spacetime derivative operator D_1
L = gen_spacetime_derivative_operator(nx, ny, nt)
# Framelet operator of the second level
W = create_framelet_operator(nx, ny, 2)
\end{lstlisting}
\end{tcolorbox}

\subsubsection{Emoji test problem}
This example considers real data of an ``emoji'' phantom measured at the University of Helsinki \cite{meaney2018tomographic}. We modify the data to determine two main limited angle problems as follows:
\begin{enumerate}
    \item \emph{Problem 1: 10 projection angles}. 
From the dataset \verb|DataDynamic_128x30.mat| we generate the problems
$\widehat\bA^{(t)}\widehat\bx^{(t)} + \widehat\be^{(t)} =\widehat\bb^{(t)}$, $t = 1,2,\dots, 33$ where $\widehat\bA^{(t)} \in \R^{2,170\times 16,384}$ %and $\widehat\bb^{(t)}\in \R^{2,170}$ 
are defined by taking 217 fan-beam projections around only 10 equidistant angles in $[0,2\pi)$. The forward operator $\bA$ has size $71,610 \times 540,672$.
\item \emph{Problem 2: 30 projection angles}. 
From the dataset \verb|DataDynamic_128x60.mat| we still generate 33 static problems, with $\widehat\bA^{(t)} \in \R^{6,510\times 16,384}$ computed from taking 217 fan-beam projections around 30 equidistant angles in $[0,2\pi)$.  Hence the dynamic forward operator $\bA$ has size $214,830 \times 540,672$.
\end{enumerate}
For both cases, the ground truth $\bx_{\rm true}$ is not available, but photographs of the shapes scanned \cite{meaney2018tomographic} are shown in Figure \ref{Fig: EmojiTrue}.

In TRIPs-Py, the main function to generate emoji data is $\verb|generate_emoji(dataset = your_dataset)|$ which allows the user to select the dataset to be either 30 or 60. Such function automatically downloads the data from the repository \footnote{\url{https://zenodo.org/records/1183532}}. If a noise level is provided as an argument in $\verb|generate_emoji()|$, then more Gaussian noise (on the top of the unknown noise already present in the recorded data) of the given level is added to the returned data. With $\verb|generate_emoji()|$ we generate data for solving both a sequence of static inverse problems of the form \eqref{eq: static}, where the operators are saved in $\rm \bA_{\rm seq}$ and the sinograms are the columns of $\bB$, and the dynamic problem of the form \eqref{eq: l2-lq}, with the blockdiagonal matrix $\bA$ and the stacked sinograms $\bb$. These are illustrated in Code 5.
   \begin{tcolorbox}[
    enhanced, attach boxed title to top center={yshift=-2mm},
    colback=cyan!10,
    colframe=black,
    colbacktitle=black,
    title = Code 5: Emoji test problem,
    text width = 15cm,
    fonttitle=\bfseries\color{white},
    boxed title style={size=small,colframe=darkspringgreen,sharp corners},
    sharp corners,]
% \begin{minted}[breaklines]{python}
\begin{lstlisting}[language=Python]
# Generate the emoji test problem with the dataset option 30
(F, d, Aseq, B, nx, ny, nt) = generate_emoji(dataset = 30)
# Generate the emoji test problem with the dataset option 60
(F, d, Aseq, B, nx, ny, nt) = generate_emoji(dataset = 60)
# generate a test problem by adding more noise to the real data
(F, d, Aseq, B, nx, ny, nt) = generate_emoji(dataset = 30, noise_level = 0.001)
\end{lstlisting}
\end{tcolorbox}

We run 100 iterations of MMGKS for solving both the static problems and the dynamic problem, with regularization parameter computed by GCV. The usage of the TRIPs-Py functions to solve the static problems is illustrated in Code 6, and the reconstructions obtained at time steps $t = 6, 14, 20, 26$ are shown in Figure \ref{Fig: emoji10_rec}, first row. The function $\verb|plot_recstructions_series()|$ can be used to display the reconstructions for both static and dynamic problems with the argument \verb|dynamic = True| or \verb|dynamic = False|, respectively. When displaying the results, we set negative solution entries to 0 as a post-processing step (note that this is different to applying nonnegativity constraints during the reconstructions as done, for instance, in \cite{buccini2020modulus, buccini2020linearized}). 
The usage of the TRIPs-Py functions to solve the dynamic problem is reported in Code 7, where we particularly focus on solvers that, within MMGKS, enforce (An)isotropic TV in space and anisotropic TV in time, as well as Group Sparsity.

More details about this test problem can be found in the demo \verb|demo_dynamic_Emoji.ipynb|.
% illustrate 
% Before finalizing this test problem, we illustrate how the user can call TV\_{\rm iso}, and GS in within MMGKS. To avoid repetition, we illustrate the methods for solving the dynamic inverse problem only on emoji test problem. MMGKS serves as a general framework in which the user can specify parameters for choosing $\rm TV_{\rm iso}$ and $\rm GS$. In particular, by setting isoTV = `isoTV' and GS = `GS' and specifying the problem dimensions by prob$\_$dims = (nx,ny, nt) we solve the dynamic problem with isotropic total variation or group sparsity. An illustration is shown in Code 8.
    \begin{tcolorbox}[
    enhanced, attach boxed title to top center={yshift=-2mm},
    colback=cyan!10,
    colframe=black,
    colbacktitle=black,
    title = Code 6: Solve and plot a series of static emoji problems,
    text width = 15cm,
    fonttitle=\bfseries\color{white},
    boxed title style={size=small,colframe=darkspringgreen,sharp corners},
    sharp corners,]
% \begin{minted}[breaklines]{python}
\begin{lstlisting}[language=Python]
# Define the regularization operator (only for space)
L = gen_first_derivative_operator_2D(nx, ny)
# For all problems, call MMGKS
for i in range(nt):
    b_vec = B[i].reshape((-1,1))
    (x_static_mmgks, info_mmgks) = MMGKS(Aseq[i], b_vec, L, pnorm=2, qnorm=1,projection_dim=1, n_iter = 100, regparam= 'gcv' , x_true=None, epsilon = 0.001)
    xx[i] = x_static_mmgks
# Plot static reconstructions
plot_recstructions_series(xx, (nx, ny, nt), dynamic = False, testproblem = 'Emoji', geome_x = 1,geome_x_small = 0,  save_imgs= False, save_path='./reconstruction/Emoji')
\end{lstlisting}
\end{tcolorbox}

\begin{figure}[ht!]
\centering
	\begin{minipage}{0.2\textwidth}
		\includegraphics[width=\textwidth]{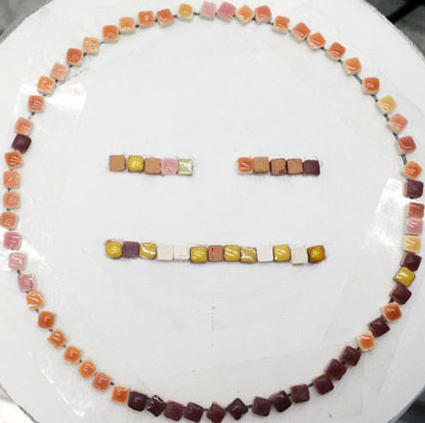}
	\end{minipage}
	\begin{minipage}{0.2\textwidth}
		\includegraphics[width=\textwidth]{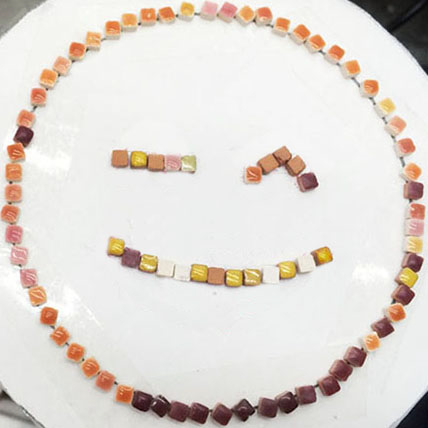}
	\end{minipage}
	\begin{minipage}{0.2\textwidth}
		\includegraphics[width=\textwidth]{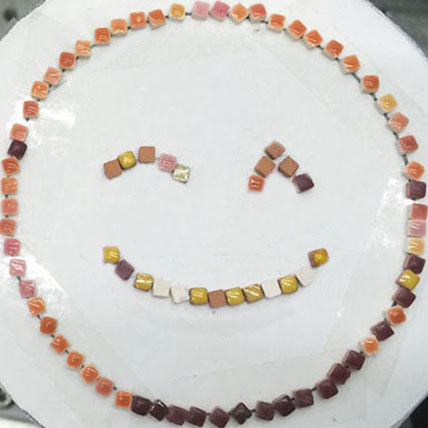}
	\end{minipage}
	\begin{minipage}{0.2\textwidth}
		\includegraphics[width=\textwidth]{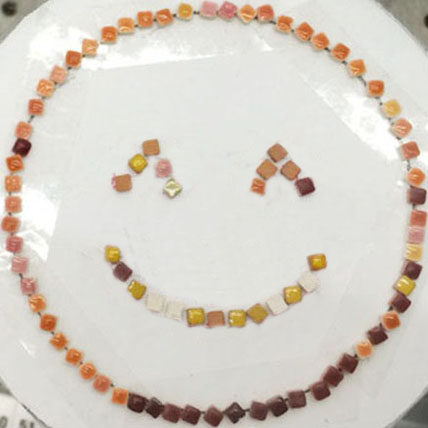}
	\end{minipage}\\
	\caption{Emoji test problem. True images \cite{meaney2018tomographic} at time instances $t = 6,14,20, 26$.}
	\label{Fig: EmojiTrue}
\end{figure}
\begin{figure}[h!]
\centering
	\begin{minipage}{0.2\textwidth}
		\includegraphics[width=\textwidth]{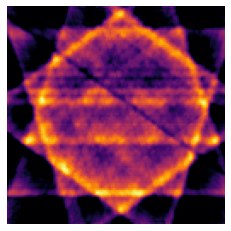}
	\end{minipage}
	\begin{minipage}{0.2\textwidth}
		\includegraphics[width=\textwidth]{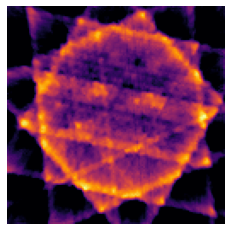}
	\end{minipage}
	\begin{minipage}{0.2\textwidth}
		\includegraphics[width=\textwidth]{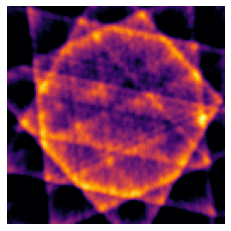}
	\end{minipage}
	\begin{minipage}{0.2\textwidth}
		\includegraphics[width=\textwidth]{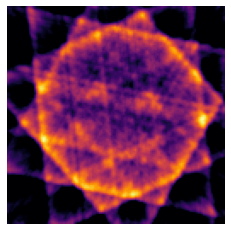}
	\end{minipage}

	\begin{minipage}{0.2\textwidth}
		\includegraphics[width=\textwidth]{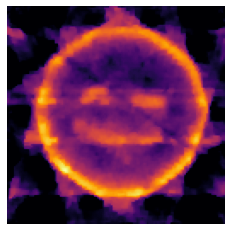}
	\end{minipage}
	\begin{minipage}{0.2\textwidth}
		\includegraphics[width=\textwidth]{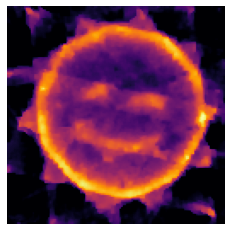}
	\end{minipage}
	\begin{minipage}{0.2\textwidth}
		\includegraphics[width=\textwidth]{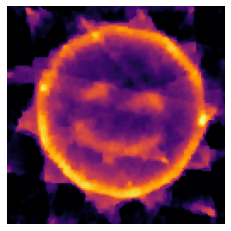}
	\end{minipage}
	\begin{minipage}{0.2\textwidth}
		\includegraphics[width=\textwidth]{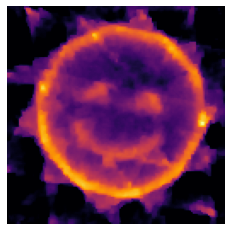}
	\end{minipage}

 	\begin{minipage}{0.2\textwidth}
		\includegraphics[width=\textwidth]{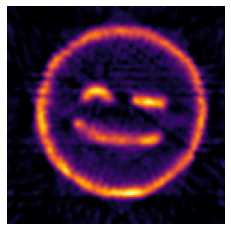}
	\end{minipage}
	\begin{minipage}{0.2\textwidth}
		\includegraphics[width=\textwidth]{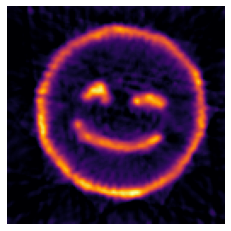}
	\end{minipage}
	\begin{minipage}{0.2\textwidth}
		\includegraphics[width=\textwidth]{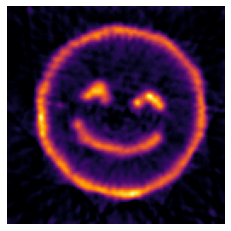}
	\end{minipage}
	\begin{minipage}{0.2\textwidth}
		\includegraphics[width=\textwidth]{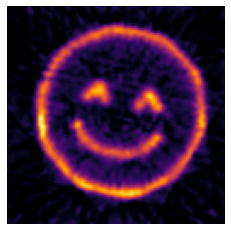}
	\end{minipage}

	\caption{Emoji test problem. Reconstruction results with $10$ projection angles. The first row shows the reconstructions with spatial anisotropic TV for the static problem, the second and third rows display the reconstructions for the dynamic problem solved with anisotropic and isotropic TV at time instances
	$t = 6, 14, 20, 26$,
	respectively, from left to right.}
	\label{Fig: emoji10_rec}
\end{figure} 

 \begin{tcolorbox}[
    enhanced, attach boxed title to top center={yshift=-2mm},
    colback=cyan!10,
    colframe=black,
    colbacktitle=black,
    title = Code 7: Solve the dynamic Emoji with different regularizers,
    text width = 15cm,
    fonttitle=\bfseries\color{white},
    boxed title style={size=small,colframe=darkspringgreen,sharp corners},
    sharp corners,]
% \begin{minted}[breaklines]{python}
\begin{lstlisting}[language=Python]
# Define the regularization operator (for space and time)
L = gen_spacetime_derivative_operator(nx, ny, nt)
data_vec = d.reshape((-1,1))
(x_mmgks, info_mmgks) = MMGKS(F, data_vec, L, pnorm=2, qnorm=1, projection_dim=1, n_iter = 100, regparam = 'gcv', x_true = None, delta = delta, epsilon = 0.001)
# Run TV_iso
(x_mmgks, info_mmgks) = MMGKS(F, data_vec, L, pnorm=2, qnorm=1, projection_dim=1, n_iter = 100, regparam = 'gcv', x_true = None, delta = delta, epsilon = 0.001, isoTV = 'isoTV', prob_dims = (nx, ny, nt)
# Run GS
(x_mmgks, info_mmgks) = MMGKS(F, data_vec, L, pnorm=2, qnorm=1, projection_dim=1, n_iter = 100, regparam = 'gcv', x_true = None, delta = delta, epsilon = 0.001, GS = 'GS', prob_dims = (nx, ny, nt)
# Plot dynamic reconstructions 
plot_recstructions_series(x_mmgks, (nx, ny, nt), dynamic = True, testproblem = 'Emoji', geome_x = 1,geome_x_small = 0,  save_imgs= False, save_path='./reconstruction/Emoji')
\end{lstlisting}
\end{tcolorbox}

%  \begin{tcolorbox}[
%     enhanced, attach boxed title to top center={yshift=-2mm},
%     colback=cyan!10,
%     colframe=black,
%     colbacktitle=black,
%     title = Code 8: Solve the dynamic problem and plot reconstructions,
%     text width = 15cm,
%     fonttitle=\bfseries\color{white},
%     boxed title style={size=small,colframe=darkspringgreen,sharp corners},
%     sharp corners,]
% \begin{minted}[breaklines]{python}
% # Solve the dynamic inverse problem
% # Define the regularization operator (for space and time)
% L = gen_spacetime_derivative_operator(nx, ny, nt)
% data_vec = d.reshape((-1,1))
% # Run TV_iso
% (x_mmgks, info_mmgks) = MMGKS(F, data_vec, L, pnorm=2, qnorm=1, projection_dim=1, n_iter = 100, regparam = 'gcv', x_true = None, delta = delta, epsilon = 0.001, isoTV = 'isoTV', prob_dims = (nx, ny, nt)
% # Run GS
% (x_mmgks, info_mmgks) = MMGKS(F, data_vec, L, pnorm=2, qnorm=1, projection_dim=1, n_iter = 100, regparam = 'gcv', x_true = None, delta = delta, epsilon = 0.001, GS = 'GS', prob_dims = (nx, ny, nt)
% # Plot dynamic reconstructions 
% plot_recstructions_series(x_mmgks, (nx, ny, nt), dynamic = True, testproblem = 'Emoji', geome_x = 1,geome_x_small = 0,  save_imgs= False, save_path='./reconstruction/Emoji')
% \end{minted}
% \end{tcolorbox}

\subsubsection{STEMPO}
Starting from the Spatio-TEmporal Motor-POwered (STEMPO) ground truth phantom from \cite{heikkila2022stempo}, within TRIPs-Py we generate both static and dynamic inverse problems, with both simulated and real data. The data for this example can be downloaded from the repository \footnote{\url{https://zenodo.org/records/8239013}}.
% \begin{enumerate}
% \item \emph{Simulated data:}
%  % we consider the true STEMPO phantoms and we compute sinograms with operators generated in ASTRA. In particular, 
 
For the simulated data we consider the STEMPO ground truth phantom \linebreak[4]
 \verb|stempo_ground_truth_2d_b4.mat|, which contains 360 images of size $560\times 560$. From this dataset we retain $n_t$ images, chosen uniformly from $1$ to $360$ with a factor of $8$, i.e., we choose the $1$st, the $8$th, ... up to the $360$th image: these represent the ground truth at $n_t$ time instances; $n_t$ can be given in input by the user. Using the ASTRA toolbox \cite{van2015astra} we generate the forward operators $\widehat\bA^{(t)}$, $t = 1,2,\dots, n_t$, each defined with respect to the angles stored in $n_t$ vectors of length $11$, with 791 parallel rays departing from each angle. Each angle vector is generated by choosing $11$ equispaced degree angles from $(5*(t-1), 5*(t-1) + 140)$, for $t = 1, 2, \dots, n_t$, which are then converted to radian. For instance, for $n_t = 20$, the forward operators at each time instant are of size ${8701 \times 313600}$, while the dynamic forward operator is a block diagonal matrix of size ${174020 \times 6272000}$. 
 % We apply the forward operators $\widehat\bA^{(t)}$, to the $n_t$ true images $\widehat\bx^{(t)}$, to obtain $n_t$ sinograms $\bd^{(t)} \in \R^{8701}$, with $\bD^{(t)} \in \R^{791 \times 11}$, for $t = 1, 2, \dots, nt$. (HOW MANY PROJECTION RAYS, WHICH GEOMETRY?)
We perturb each measured vectorized sinogram $\widehat\bd_{\rm true}^{(t)}$ with white Gaussian noise of level given in input by the user. , i.e.,  the noise vector $\be^{(t)}$ has mean zero and a rescaled identity covariance matrix. 
% We refer to the ratio 
% $\sigma^{(t)} = \|\be^{(t)}\|_{2}/\|\bA^{(t)}\bx^{(t)}\|_{2}$ as the noise level -- THIS COULD PROBABLY BE MORE GENERAL (MAYBE FOR THE ADDNOISE FCN, TO BE DESCRIBED SOMEWHERE?). 
The true images $\bx^{(t)}$ at time steps $n_t = 1, 10, 20, 30$ are shown in Figure ~\ref{Figure: TrueStempo}. 
% \item \emph{Real data:} For the real data scenario the user can select data provided from \cite{heikkila2022stempo}.
% \end{enumerate}
An illustration for both real and simulated gata generation is shown in Code 9. More details about this test problem can be found in the demo \verb|demo_dynamic_Stempo.ipynb|.
% \end{enumerate}
   \begin{tcolorbox}[
    enhanced, attach boxed title to top center={yshift = -2mm},
    colback = cyan!10,
    colframe=black,
    colbacktitle=black,
    title = Code 9: Stempo test problem,
    text width = 15cm,
fonttitle=\bfseries\color{white},
    boxed title style={size=small,colframe=darkspringgreen,sharp corners},
    sharp corners,]
% \begin{minted}[breaklines]{python}
\begin{lstlisting}[language=Python]
# Generate a simulated STEMPO test problem
(F, d, Aseq, B, nx, ny, nt, savedelta, truth) = generate_stempo(data_set = 'simulation', data_thinning = 2, nt = 10)
# Generate a STEMPO test problem with real data
(F, d, Aseq, B, nx, ny, nt, savedelta, truth) = generate_stempo(data_set = 'real', data_thinning = 2)
\end{lstlisting}
\end{tcolorbox}

% \begin{figure}[!ht]
% \centering
% {\includegraphics[width = 0.2\textwidth]{Figures/STEMPO/stempo_true_1.png}}
% \includegraphics[width = 0.2\textwidth]{Figures/STEMPO/stempo_true_10.png} 
% \includegraphics[width = 0.2\textwidth]{Figures/STEMPO/stempo_true_15.png}
% \includegraphics[width = 0.2\textwidth]{Figures/STEMPO/stempo_true_19.png}
% \caption{STEMPO test problem: True images at time steps $t = 2, 11, 16, 20$.}
% \label{Figure: TrueStempo}
% \end{figure}

\begin{figure}[h!]
\centering
	\centering
\begin{tabular}{cccc}
		\includegraphics[width=0.17\textwidth]{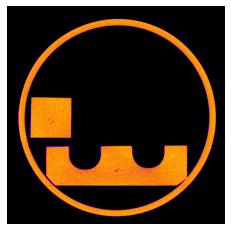}&
		\includegraphics[width=0.17\textwidth]{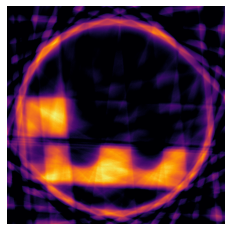}&
		\includegraphics[width=0.17\textwidth]{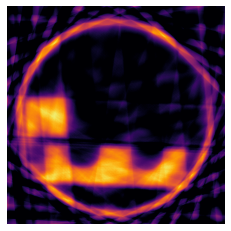}&
		\includegraphics[width=0.17\textwidth]{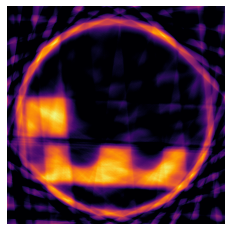}\\
(a) & (b)  & (c)  & (d) 
\end{tabular}
	\caption{STEMPO test problem. (a) True image, and reconstructions at time instance $t=10$ by (b) static anisotropic TV, (c) Hybrid LSQR, (d) dynamic isotropic TV.}
\label{Figure: TrueStempo}
\end{figure}

\subsubsection{Cross Phantom}
The last example considers real data of the cross phantom measures at the University of Helsinki \cite{latva2018tomographic}. We obtain the forward operator and the data from the file \verb|DataDynamic_128x15.mat|, which can be downloaded from the repository \footnote{\url{https://zenodo.org/records/1446516}}. The spatial resolution is $128\times 128$ pixels, while the time resolution is 16. The CT sinogram that represents the available data consists of $16$ time frames that are generated by measuring a 2D cross-section of the dynamic phantom that is built from an aluminum stick, a graphite stick, and candle wax. The measurements were collected by considering 15 projection angles at each time instance, with angles shifted by one degree from a given time instance to the next one, i.e., if the angles for time step $t_1$ are $[1,15,29,..]$, then for time step $t_2$ they are $[2,16,30,..]$. The measurement matrix that represents the forward operator obtained with a cone-beam geometry is a sparse matrix of size $33600 \times 262 144$, and the sinogram matrix is represented as $\bB \in \R^{140 \times 240}$, i.e., 15 projections on 16 images. The instance of the Cross Phantom test problem just described can be easily generated in TRIPs-Py as shown in Code 1Reconstructed images with MMGKS for the static inverse problems at time instances $t = 1, 5, 10, 15$ are shown in the first row of Figure \ref{Figure: CrossPhantom}. The second row of Figure \ref{Figure: CrossPhantom} shows reconstructed images by solving the dynamic inverse problem with anisotropic TV at time instances  $t = 1, 5, 10, 15$ from left to right, respectively. More details about this test problem can be found in the demo \verb|demo_dynamic_CrossPhantom.ipynb|. 

   \begin{tcolorbox}[
    enhanced, attach boxed title to top center={yshift=-2mm},
    colback=cyan!10,
    colframe=black,
    colbacktitle=black,
    title = Code 10: CrossPhantom test problem,
    text width = 15cm,
fonttitle=\bfseries\color{white},
    boxed title style={size=small,colframe=darkspringgreen,sharp corners},
    sharp corners,]
% \begin{minted}[breaklines]{python}
\begin{lstlisting}[language=Python]
# Generate the CrossPhantom test problem
(F, d, Aseq, B, nx, ny, nt) = generate_crossPhantom(dataset = 15)
\end{lstlisting}
\end{tcolorbox}

\begin{figure}[!ht]
\centering
{\includegraphics[width = 0.24\textwidth,angle =-90]{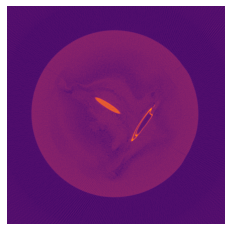}}
{\includegraphics[width = 0.24\textwidth,angle =-90]{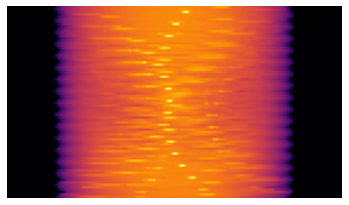}}
\caption{CrossPhantom test problem. Left frame: the high-resolution  filtered back-projection reconstruction of the time-dependent cross phantom computed from 360 projections. Right frame: sinogram that contains measurements from $80$ angels for 16 time instances together.}
	\label{Figure: CrossPhantom}
\end{figure}

\begin{figure}[h!]
\centering
	\begin{minipage}{0.2\textwidth}
		\includegraphics[width=\textwidth]{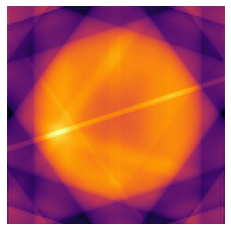}
	\end{minipage}
	\begin{minipage}{0.2\textwidth}
		\includegraphics[width=\textwidth]{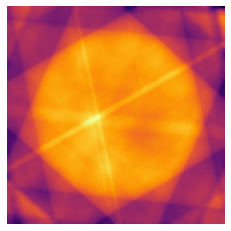}
	\end{minipage}
	\begin{minipage}{0.2\textwidth}
		\includegraphics[width=\textwidth]{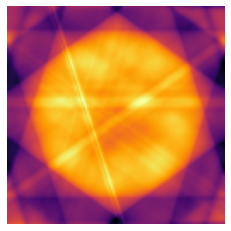}
	\end{minipage}
	\begin{minipage}{0.2\textwidth}
		\includegraphics[width=\textwidth]{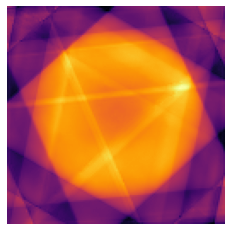}
	\end{minipage}

	\begin{minipage}{0.2\textwidth}
		\includegraphics[width=\textwidth]{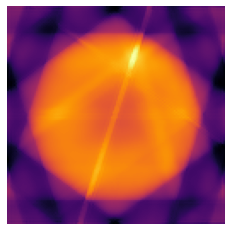}
	\end{minipage}
	\begin{minipage}{0.2\textwidth}
		\includegraphics[width=\textwidth]{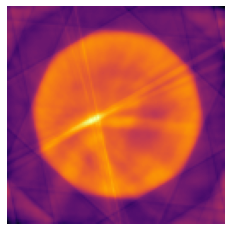}
	\end{minipage}
	\begin{minipage}{0.2\textwidth}
		\includegraphics[width=\textwidth,]{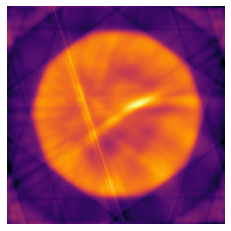}
	\end{minipage}
	\begin{minipage}{0.2\textwidth}
		\includegraphics[width=\textwidth]{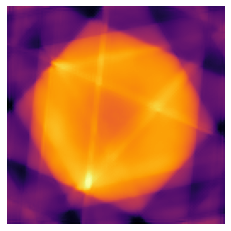}
	\end{minipage}
	\caption{CrossPhantom test problem. Reconstructions with MMGKS for the static test problems (first row) and anisotropic TV for the dynamic inverse problems (second row) at time steps $t=1, 5, 10, 15$, from left to right. }
	\label{Fig: crossPhantom_reconstructions}
\end{figure}

\section{Conclusions and outlook}\label{sec: conclusions}
In this paper we describe TRIPs-Py, a Python package collecting solvers for small and large-scale discrete ill-posed inverse problems, and test problems. The package allows the user to easily set-up a test problem, with both simulated and real data, and employ built-in regularization methods to solve the inverse problem. Among solvers for small-scale problems, TRIPs-Py includes direct regularization methods such as truncated (G)SVD and Tikhonov. Among iterative regularization, TRIPs-Py includes methods based on both standard and generalized Krylov subspace methods. 
% GMRES, Arnoldi Tikhonov, Golub Kahan Tikhonov, hybrid methods (Hybid\_GMRES and Hybrid\_LSQR) that rely on Krylov subspaces and other solvers such as MMGKS, \rm \textit{}TV\_{\rm aniso}, and TV\_{GS} that rely on generalized Krylov subspaces. 
% We provide a set of simulated and real data problems. 
Among the test problems, special emphasis is given to a framework for handling dynamic simulated and real data. We envision that the user can provide their own data and use our solution methods within TRIPs-Py, as well as test their solution methods with the TRIPs-Py test problems. 

This is a first public version of TRIPs-Py. Future developments will focus on adding new solution methods, functionalities, and improve existing data generation and visualisation tools. Specifically, in future versions of TRIPs-Py, we plan to:
\begin{enumerate}
% \item [\diamond] Enrich the package with other Krylov subspace and regularization methods.
\item [$\diamond$] Include solvers based on flexible Krylov subspace methods to, e.g., enforce sparsity on the desired solution; see \cite{gazzola2021iteratively}.
\item [$\diamond$]  Include methods that can enforce constraints into the solutions; see, e.g., \cite{NNFCGLS, buccini2020modulus, buccini2020linearized}.
\item [$\diamond$]  Implement memory-aware iterative methods for edge-preserving and sparsity for solving large-scale and massive inverse problems; see \cite{pasha2023recycling}.
\item [$\diamond$]  Add other strategies to automatically determine the regularization parameter(s).
\item [$\diamond$] Add other instances of the existing deblurring test problems (e.g., a variety of PSFs) and tomography test problems (e.g., different scanning geometries). 
\item [$\diamond$]  Add other test problems, such as hyperspectral imaging and spectral CT.
\item [$\diamond$] Develop tools for high-order representation of dynamic imaging data.
\end{enumerate}

Users are welcomed to contribute to TRIPs-Py by solution methods, functionalities, and/or test problems and data.

\section*{Acknowledgments}
MP gratefully acknowledges support from the NSF under award No.\ 2202846. MP further acknowledges partial support from the NSF-AWM Mentoring Travel award.
Both MP and SG acknowledge the Isaac Newton Institute for Mathematical Sciences, Cambridge, for the support and hospitality during the programme ``Rich and Nonlinear Tomography - a multidisciplinary approach" (supported by EPSRC grant no EP/R014604/) where partial work on this paper was undertaken. SG would like to thank Ludovico Carozza for his advice on many Python functionalities. We would like to thank Tatiana Bubba for discussions about avoiding `inverse crimes' in inverse problems.

\bibliographystyle{vancouver}
\bibliography{TRIPS_ARXIV_January30.bib}

\begin{thebibliography}{10}

\bibitem{boas2001imaging}
Boas DA, Brooks DH, Miller EL, DiMarzio CA, Kilmer M, Gaudette RJ, et~al.
\newblock Imaging the body with diffuse optical tomography.
\newblock IEEE signal processing magazine. 2001;18(6):57-75.

\bibitem{miller2012environmental}
Miller EL, Abriola LM, Aghasi A.
\newblock Environmental remediation and restoration: hydrological and geophysical processing methods.
\newblock IEEE Signal Processing Magazine. 2012;29(4):16-26.

\bibitem{bennett1996generalized}
Bennett AF, Chua BS, Leslie L.
\newblock Generalized inversion of a global numerical weather prediction model.
\newblock Meteorology and Atmospheric Physics. 1996;60(1):165-78.

\bibitem{hansenRT}
Hansen PC.
\newblock Regularization tools: A Matlab package for analysis and solution of discrete ill-posed problems.
\newblock Numer Algo. 1994;6(3):1-35.

\bibitem{gazzola2019ir}
Gazzola S, Hansen PC, Nagy JG.
\newblock IR Tools: a MATLAB package of iterative regularization methods and large-scale test problems.
\newblock Numerical Algorithms. 2019;81(3):773-811.

\bibitem{NaPaPe04}
Nagy J, Palmer K, Perrone L.
\newblock Iterative Methods for Image Deblurring: A {MATLAB} Object Oriented Approach.
\newblock Numerical Algorithms. 2004;36(1):73-93.

\bibitem{hansenAIRII}
Hansen PC, J\o{}rgensen JS.
\newblock AIR Tools II: Algebraic iterative reconstruction methods, improved implementation.
\newblock Numer Algo. 2018;79:107-37.

\bibitem{TIGRE}
Biguri A, Dosanjh M, Hancock S, Soleimani M.
\newblock {TIGRE}: a {MATLAB-GPU} toolbox for {CBCT} image reconstruction.
\newblock Biomed Phys Eng Express. 2016;2:055010.

\bibitem{CILI}
J\o{}rgensen JS, Ametova E, Burca G, Fardell G, Papoutsellis E, Pasca E, et~al.
\newblock {C}ore {I}maging {L}ibrary - Part {I}: a versatile {P}ython framework for tomographic imaging.
\newblock Phil Trans R Soc. 2021;A:3792020019220200192.

\bibitem{van2015astra}
Van~Aarle W, Palenstijn WJ, De~Beenhouwer J, Altantzis T, Bals S, Batenburg KJ, et~al.
\newblock The ASTRA Toolbox: A platform for advanced algorithm development in electron tomography.
\newblock Ultramicroscopy. 2015;157:35-47.

\bibitem{ODL}
et~al AJ. odlgroup/odl: ODL 0.7.0. Zenodo; 2018.
\newblock \url{doi:10.5281/zenodo.1442734}.

\bibitem{hansen2010discrete}
Hansen PC.
\newblock Discrete inverse problems: Insight and algorithms.
\newblock SIAM; 2010.

\bibitem{hansen1998rank}
Hansen PC.
\newblock Rank-deficient and discrete ill-posed problems: numerical aspects of linear inversion.
\newblock SIAM; 1998.

\bibitem{saad2003iterative}
Saad Y.
\newblock Iterative methods for sparse linear systems.
\newblock SIAM; 2003.

\bibitem{review2023}
Chung J, Gazzola S.
\newblock Computational methods for large-scale inverse problems: a survey on hybrid projection methods.
\newblock arXiv preprint arXiv:210507221. 2023.

\bibitem{fenu2016gcv}
Fenu C, Reichel L, Rodriguez G.
\newblock GCV for Tikhonov regularization via global Golub--Kahan decomposition.
\newblock Numerical Linear Algebra with Applications. 2016;23(3):467-84.

\bibitem{lewis2009arnoldi}
Lewis B, Reichel L.
\newblock Arnoldi--Tikhonov regularization methods.
\newblock Journal of Computational and Applied Mathematics. 2009;226(1):92-102.

\bibitem{lampe2012large}
Lampe J, Reichel L, Voss H.
\newblock Large-scale Tikhonov regularization via reduction by orthogonal projection.
\newblock Linear algebra and its applications. 2012;436(8):2845-65.

\bibitem{lanza2015generalized}
Lanza A, Morigi S, Reichel L, Sgallari F.
\newblock A generalized {K}rylov subspace method for $\ell_p-\ell_q$ minimization.
\newblock SIAM Journal on Scientific Computing. 2015;37(5):S30-50.

\bibitem{pasha2021efficient}
Pasha M, Saibaba AK, Gazzola S, Espanol MI, de~Sturler E.
\newblock Efficient edge-preserving methods for dynamic inverse problems.
\newblock arXiv preprint arXiv:210705727. 2021.

\bibitem{GVL12}
Golub GH, Van~Loan CF.
\newblock Matrix {C}omputations, 4th ed.
\newblock Johns Hopkins University Press, Baltimore; 2013.

\bibitem{elden1982weighted}
Eld{\'e}n L.
\newblock A weighted pseudoinverse, generalized singular values, and constrained least squares problems.
\newblock BIT Numerical Mathematics. 1982;22(4):487-502.

\bibitem{hansen2006deblurring}
Hansen PC, Nagy JG, O'leary DP.
\newblock Deblurring images: matrices, spectra, and filtering.
\newblock SIAM; 2006.

\bibitem{björck2014numerical}
Bj{\"o}rck {\AA}.
\newblock Numerical Methods in Matrix Computations.
\newblock Texts in Applied Mathematics. Springer International Publishing; 2014.

\bibitem{arnoldi1951principle}
Arnoldi WE.
\newblock The principle of minimized iterations in the solution of the matrix eigenvalue problem.
\newblock Quarterly of applied mathematics. 1951;9(1):17-29.

\bibitem{lanczos1950iteration}
Lanczos C.
\newblock An iteration method for the solution of the eigenvalue problem of linear differential and integral operators.
\newblock United States Governm. Press Office Los Angeles, CA; 1950.

\bibitem{lange2016mm}
Lange K.
\newblock MM optimization algorithms.
\newblock SIAM; 2016.

\bibitem{rodriguez2008efficient}
Rodriguez P, Wohlberg B.
\newblock An efficient algorithm for sparse representations with $\ell_p$ data fidelity term.
\newblock In: Proceedings of 4th IEEE Andean Technical Conference (ANDESCON); 2008. .

\bibitem{huang2017majorization}
Huang G, Lanza A, Morigi S, Reichel L, Sgallari F.
\newblock Majorization--minimization generalized Krylov subspace methods for $\ell_p-\ell_q$ optimization applied to image restoration.
\newblock BIT Numerical Mathematics. 2017;57(2):351-78.

\bibitem{bach2012optimization}
Bach F, Jenatton R, Mairal J, Obozinski G.
\newblock Optimization with Sparsity-Inducing Penalties.
\newblock Foundations and Trends in Machine Learning. 2012;4(1):1-106.

\bibitem{ReichelAnzf}
Reichel L, Shyshkov A.
\newblock A new zero-finder for Tikhonov regularization.
\newblock BIT Numerical Mathematics. 2008;48(3):627-43.

\bibitem{gazzola2020survey}
Gazzola S, Landman MS.
\newblock Krylov methods for inverse problems: Surveying classical, and introducing new, algorithmic approaches.
\newblock Mitteilungen der Gesellschaft für Angewandte Mathematik und Mechanik. 2020;43(4).

\bibitem{BRXX}
Buccini A, Reichel L.
\newblock An $\ell^2$-$\ell^q$ Regularization Method for Large Discrete Ill-Posed Problems.
\newblock Journal of Scientific Computing. 2019;78:1526-49.

\bibitem{NoRu14}
Novati P, Russo MR.
\newblock A {GCV}-based {A}rnoldi-{T}ikhonov regularization method.
\newblock BIT. 2014;54:501-21.

\bibitem{bucciniGCV}
Buccini A, Reichel L.
\newblock Generalized cross validation for $\ell_p-\ell_q$ minimization.
\newblock Numer Algor. 2021;88:1595–1616.

\bibitem{buccini2020modulus}
Buccini A, Pasha M, Reichel L.
\newblock Modulus-based iterative methods for constrained $\ell_p-\ell_q$ minimization.
\newblock Inverse Problems. 2020;36(8):084001.

\bibitem{buccini2020linearized}
Buccini A, Pasha M, Reichel L.
\newblock Linearized {K}rylov subspace {B}regman iteration with nonnegativity constraint.
\newblock Numerical Algorithms. 2020:1-24.

\bibitem{COS09b}
Cai JF, Osher S, Shen Z.
\newblock Linearized {Bregman} iterations for frame-based image deblurring.
\newblock SIAM Journal on Imaging Sciences. 2009;2:226-52.

\bibitem{mueller2012linear}
Mueller JL, Siltanen S.
\newblock Linear and nonlinear inverse problems with practical applications.
\newblock SIAM; 2012.

\bibitem{PCHCT}
Hansen PC, Jørgensen JS, Lionheart WRB.
\newblock Computed Tomography: Algorithms, Insight, and Just Enough Theory.
\newblock SIAM, Philadelphia; 2021.

\bibitem{lan2023spatiotemporal}
Lan S, Pasha M, Li S.
\newblock Spatiotemporal {B}esov Priors for {B}ayesian Inverse Problems.
\newblock arXiv preprint arXiv:230616378. 2023.

\bibitem{meaney2018tomographic}
Meaney A, Purisha Z, Siltanen S.
\newblock Tomographic X-ray data of 3D emoji.
\newblock arXiv preprint arXiv:180209397. 2018.

\bibitem{heikkila2022stempo}
Heikkil{\"a} T.
\newblock STEMPO--dynamic X-ray tomography phantom.
\newblock arXiv preprint arXiv:220912471. 2022.

\bibitem{latva2018tomographic}
Latva-{\"A}ij{\"o} S, Meaney A, Siltanen S.
\newblock Tomographic X-ray data of 3D cross phantom.
\newblock arXiv preprint arXiv:180900166. 2018.

\bibitem{gazzola2021iteratively}
Gazzola S, Nagy JG, Sabaté~Landman M.
\newblock Iteratively Reweighted FGMRES and FLSQR for Sparse Reconstruction.
\newblock SIAM Journal on Scientific Computing. 2021;(0):S47-69.

\bibitem{NNFCGLS}
Gazzola S, Wiaux Y.
\newblock Fast nonnegative least squares through flexible {K}rylov subspaces.
\newblock SIAM J Sci Comput. 2017;39:A655-79.

\bibitem{pasha2023recycling}
Pasha M, de~Sturler E, Kilmer ME.
\newblock Recycling MMGKS for large-scale dynamic and streaming data.
\newblock arXiv preprint arXiv:230915759. 2023.

\end{thebibliography}
\end{document}